\documentclass[10pt,twocolumn,twoside]{IEEEtran}

\IEEEoverridecommandlockouts        

\usepackage{etex}
\usepackage[dvipsnames]{xcolor}
\usepackage[vlined,ruled]{algorithm2e}
\usepackage{cite}
\usepackage{booktabs}

\usepackage{graphicx,color}
\graphicspath{{./fig/}}
\usepackage{epstopdf}

\usepackage{amsmath}
\usepackage{amssymb}
\usepackage{mathrsfs}
\usepackage{dsfont}

\usepackage[font=footnotesize]{caption}
\usepackage{subcaption}

\newtheorem{theorem}{Theorem}[section]
\newtheorem{lemma}[theorem]{Lemma}
\newtheorem{definition}{Definition}
\newtheorem{corollary}[theorem]{Corollary}

\newtheorem{remark}{Remark}

\newcommand*\fvec[1]{\ensuremath{\mathbf{#1}}}                           
\newcommand*\mc[0]{\mathcal}                                                  	
\DeclareMathOperator*{\minimize}{minimize} 									
\newcommand{\transpose}{\mathsf{T}} 
\newcommand{\Span}{\mathrm{Span\,}}											
\newcommand{\subscr}[2]{#1_{\textup{#2}}}										
\newcommand{\supscr}[2]{#1^{\textup{#2}}}										
\newcommand{\map}[3]{#1: #2 \rightarrow #3}									

\DeclareMathOperator*{\diag}{\mathrm{diag}}									    
\DeclareMathOperator*{\sinbf}{\mathrm{\bf sin}} 	                           

\newcommand{\mycircle}{\ensuremath{\mathbb S^{1}}}

\newcommand{\torus}{\ensuremath{\mathbb T}}
\newcommand{\real}{\mathbb{R}}
\newcommand{\complex}{\mathbb{C}}
\newcommand{\rot}{\operatorname{rot}}

\newcommand{\until}[1]{\{1,\dots, #1\}}												

\DeclareMathAlphabet{\mathpzc}{OT1}{pzc}{m}{it}

\newcommand\oprocendsymbol{\hbox{$\square$}}
\newcommand\oprocend{\relax\ifmmode\else\unskip\hfill\fi\oprocendsymbol}

\def\QEDopen{{\setlength{\fboxsep}{0pt}\setlength{\fboxrule}{0.2pt}\fbox{\rule[0pt]{0pt}{1.3ex}\rule[0pt]{1.3ex}{0pt}}}}
\def\QED{\QEDopen} 

\renewcommand{\theenumi}{(\roman{enumi})}

\title{Breaking the Hierarchy: Distributed Control \\\& Economic Optimality in 
Microgrids %
\thanks{
This work was supported in part by ETH startup funds, the National Science and Engineering Research Council of Canada, and the National Science Foundation NSF CNS-1135819. A preliminary version of part of the results in this document has been presented in \cite{HB-JWSP-FD-FB:13j}.}}

\author{Florian D\"orfler, John W. Simpson-Porco, and Francesco Bullo%
    \thanks{%
     Florian D{\"o}rfler is with the Automatic Control Laboratory, Swiss Federal Institute of Technology (ETH) Z\"urich, Switzerland. Email: {\tt dorfler@ethz.ch}. %
    J. W. Simpson-Porco and F. Bullo are with the Mechanical Engineering Department, University~of~California~Santa Barbara. Email: {\tt \{johnwsimpsonporco,bullo\}@engineering.ucsb.edu}. %
}}

% ...

\begin{document}
\maketitle
\thispagestyle{empty}
\pagestyle{empty}


\begin{abstract}
Modeled after the hierarchical control architecture of power transmission systems, a layering of primary, secondary, and tertiary control has become the standard operation paradigm for {islanded} microgrids. Despite this superficial similarity, the control objectives in microgrids across these three layers are varied and ambitious, 
and they must be achieved while allowing for  robust plug-and-play 
operation and maximal flexibility, without hierarchical decision making and
time-scale separations. In this work, we explore control strategies for
these three layers and illuminate some possibly-unexpected connections and
dependencies among them. Building from a first-principle analysis of
decentralized primary droop control, we study centralized, decentralized,
and distributed architectures for secondary frequency regulation. We find
that averaging-based distributed controllers using communication among the
generation units offer the best combination of flexibility and
performance. We further leverage these results to study constrained AC
economic dispatch in a tertiary control layer. Surprisingly, we show that
the minimizers of the economic dispatch problem are in one-to-one
correspondence with the set of steady-states reachable by 
droop control. In other words, the adoption of droop control is necessary
and sufficient to achieve economic optimization. 
This equivalence results in simple guidelines to select the droop coefficients, which include the known criteria for power sharing. 
We illustrate the performance and robustness of our designs through simulations.%
\end{abstract}


\section{Introduction}
\label{Section: Introduction}

With the goal of integrating distributed renewable generation and energy storage systems, the concept of a {\em microgrid} has recently gained popularity \cite{RHL:02,MCC-DMD-RA:93,QCZ-TH:13,JMG-JCV-JM-LGDV-MC:11}. Microgrids are low-voltage electrical distribution networks, heterogeneously composed of distributed generation, storage, load, and managed autonomously from the larger transmission network. Microgrids are able to connect to a larger power system, 
but are also able to island themselves and operate independently. {Such islanding could be the result of an emergency, such as an outage of the larger utility grid, or may be by design in an isolated grid.}

Distributed energy sources in a microgrid generate either DC  or variable frequency AC power, and are interfaced with an AC grid via power electronic DC/AC \emph{inverters}. {In islanded operation,} it is through these inverters, cooperative actions must be taken to ensure frequency synchronization, voltage stability, power balance, load sharing, regulation of disturbances, {and economic operation}~\cite{JMG-MC-TLL-PLC:13,JAPL-CLM-AGM:06}. A variety of control and decision architectures \textemdash{} ranging from centralized to fully decentralized \textemdash{} have been proposed to address these challenges \cite{AM-EO-DM-OO:10,JAPL-CLM-AGM:06,JMG-JCV-JM-LGDV-MC:11,JMG-MC-TLL-PLC:13}. 
 In transmission networks, 
the different control tasks are separated in their time scales and aggregated into a hierarchy. Similar operation layers have been proposed for microgrids.

\textit{Control Hierarchy in Transmission Systems:} %
 The foundation of this hierarchy, termed {\em primary control}, must rapidly balance generation and demand, while sharing the load, synchronizing the AC voltage frequencies, and stabilizing their magnitudes. This is accomplished via decentralized {\em droop} control, where generators are controlled such that their power injections are proportional to their voltage frequencies and magnitudes \cite{PK:94}. 

Droop controllers induce steady-state errors in frequency and voltage magnitudes, which  are corrected in a {\em secondary control} layer. At the transmission level, the network is partitioned into control areas, and a few selected generators then balance local generation in each area with load and inter-area power transfers. Termed {\em automatic generation control} (AGC), 
this architecture is based on centralized integral control and operates on a slower time scale than primary control \cite{JM-JWB-JRB:08}.

The operating point stabilized by primary and secondary control is scheduled in a \emph{tertiary} control layer, to establish fair load sharing among the sources, or to dispatch the generation to minimize operational costs. In conventional operation of bulk power systems, an economic dispatch is optimized offline, in a centralized fashion, using precise load forecasts \cite{AJW-BFW:96}. In {\cite{MA-DVD-KHJ-HS:13,MA-DVD-HS-KHJ:14,EM-SHL:14,SY-LC:14,CZ-UT-NL-SL:13,NL-CL-ZC-SHL:13,XZ-AP:13,MB-CDP-ST:14}} it has been shown that the dynamics of a power transmission system with synchronous generators and AGC naturally optimize variations of the economic dispatch. 

\textit{Adaption of Control Layers to Microgrids:}
{With regards to primary control in islanded microgrids, inverters are typically controlled to emulate the droop characteristics of synchronous generators \cite{MCC-DMD-RA:93,JMG-JCV-JM-LGDV-MC:11,JAPL-CLM-AGM:06,JMG-MC-TLL-PLC:13,QCZ-TH:13}.
Despite forming the foundation of microgrid operation, networks of droop-controlled inverters have only recently been subject to a rigorous analysis \cite{JWSP-FD-FB:12u,JWSP-FD-FB:13h}. We also refer to \cite{ZW-MX-MDL:13,VM-FV:13,NA-SG:13,LYL:13,JS-RO-AA-JR-TS:13-updated,JS-DG-JR-TS:13,UM-MM:14} for further results.} 
{To compensate for steady-state deviations induced by droop control, secondary integral control strategies akin to AGC have been adapted to microgrids. 
Whereas fully decentralized integral controllers successfully regulate the frequency, they result in steady-state deviations from the desired power\,\,injection profile \cite{NA-SG:13b}. 
Thus, distributed controllers merging primary and secondary control have been proposed based upon continuous-time averaging with all-to-all \cite{QS-JG-JMV:13-updated,HL-BJC-WZ-XS:13} or nearest-neighbor \cite{LYL:13,JWSP-FD-FB:12u} communication. 
In transmission grids, the tertiary optimization layer can be merged with the primary and secondary layer based on continuous-time optimization approaches \cite{MA-DVD-KHJ-HS:13,MA-DVD-HS-KHJ:14,EM-SHL:14,SY-LC:14,XZ-AP:13,CZ-UT-NL-SL:13,NL-CL-ZC-SHL:13,MB-CDP-ST:14}. 
Similar discrete-time approaches are based on game-theoretic ideas \cite{MNE-CAM-NQ:14} or discrete-time averaging\,\,algorithms \cite{RM-SD-BBC:12,STC-ADGC-CNH:13}.
}

{
The above approaches employ varying models ranging from linear to nonlinear differential-algebraic, some of which are not appropriate in microgrids such as lossless lines and rotational inertia. Some of the proposed strategies are validated only numerically without providing further analytic insights. Often the primary and secondary control loops may interact in an adverse way unless a time-scale separation is enforced, the gains are carefully tuned, or an estimate of the load is known.%
}

\textit{Transmission Level vs. Distribution \& Microgrids: }
While the hierarchical architecture has been adapted from the transmission level to microgrids, the control challenges and architecture limitations imposed by the microgrid framework are as diverse as they are daunting. The low levels of inertia in microgrids mean that primary control must be fast\,and\,reliable to maintain voltages, frequencies, and power flows within acceptable tolerances, while the highly variable and distributed nature of microgrids preclude centralized control strategies\,of any kind. Microgrid controllers must be able to adapt in real time to unknown and variable loads and network conditions. In short, the three layers of the control hierarchy must allow for as close to plug-and-play operation as possible, be either distributed or completely decentralized, {without knowledge\,of the system model and the load and generation\,profile}, and operate seamlessly without a pre-imposed separation of time\,scales.

\textit{Contributions and Contents: }
%
In Section \ref{Section: Microgrids and their Control Challenges}, we present\,\,a comprehensive modeling {and control} framework for
microgrids with heterogenous components and control objectives. 
{Our approach builds on a first-principle nonlinear differential-algebraic model, decentralized primary droop controllers, and networks with constant (not necessarily zero) resistance-to-reactance ratios extending the conventional lossless models.

In Section \ref{Section: Primary Control Strategies}, we review the properties and limitations of  droop control, including the conditions for the existence of stable synchronized solutions satisfying actuation constraints and proportional load sharing.}
{Moreover, we show the following reachability result: the set of feasible setpoints
for generation dispatch} is in one-to-one correspondence with the
set of steady-states reachable via decentralized droop control.

In Section \ref{Section: Secondary Control Strategies}, we study several decentralized and distributed secondary integral control strategies. We first discuss the limitations of decentralized secondary integral control akin to AGC. Next, we study distributed secondary control strategies {based on averaging}. We provide a rigorous analysis for the strategies proposed in \cite{QS-JG-JMV:13-updated,HL-BJC-WZ-XS:13} for a proper choice of control gains and compare them to our earlier work \cite{JWSP-FD-FB:12u} with regards to tuning limitations and communication complexity. We show that all these distributed strategies successfully regulate the frequency, maintain the injections and stability properties of the primary droop controller, and do not require any separation of time scales. Finally, we demonstrate that these properties are maintained when only a subset of generating units participate in secondary control. The effectiveness {and practical applicability} of the proposed distributed secondary control strategies {has been validated experimentally; see \cite{JWSP-FD-QS-JMG-FB:13e,QS-JG-JMV:13-updated,JWSP-QS-FD-JMV-JMG-FB:14}}.

In Section \ref{Section: Tertiary Control Strategies}, we study tertiary
control policies that minimize an economic dispatch problem. We leverage a
recently discovered relation between AC and DC power flows
\cite{FD-MC-FB:11v-pnas,FD-FB:13c} and show that the set of minimizers of
the nonlinear and non-convex AC economic dispatch optimization problem are
in one-to-one correspondence with the minimizers of a convex DC dispatch
problem. Our next result shows a surprising symbiotic relationship between
primary/secondary control and tertiary. We show that the minimum of the AC economic dispatch can be
achieved by a decentralized droop control design. {Whereas similar conditions are known for related transmission system problems \cite{MA-DVD-KHJ-HS:13,MA-DVD-HS-KHJ:14,EM-SHL:14,SY-LC:14,XZ-AP:13,CZ-UT-NL-SL:13,NL-CL-ZC-SHL:13,MB-CDP-ST:14} (in a simplified linear and convex setting with lossless DC power flows), we also establish a converse result:} every droop
controller results in\,a steady-state which is the minimizer of some AC
economic dispatch. We deduce, among others, that the optimal droop coefficients are inversely proportional to the marginal generation costs{, and the conventional power sharing objectives is a particular case.} 

{
In summary, we demonstrate that simple distributed and averaging-based PI controllers are able to simultaneously addresses primary, secondary, and tertiary-level objectives, without time-scale separation, relying only on local measurements and nearest-neighbor communication, and in a model-free fashion independent of the network parameters and topology, the loading profile, and the number of sources. Thus, our control strategy is suited for true plug-and-play implementation.}

In Section \ref{Section: Simulation Scenario}, we illustrate the performance and robustness of our controllers with a simulation study of the IEEE 37 bus distribution network. Finally, Section \ref{Section: Conclusions} concludes the paper. The remainder of this section introduces some preliminaries. 

\subsection*{Preliminaries and Notation}

\emph{Vectors and matrices}: Given a finite set $\mathcal{V}$, let $|\mathcal{V}|$ denote its cardinality. Given a finite index set $\mathcal{I}$ and a real-valued one-dimensional array $\{x_1,\ldots,x_{|\mathcal{I}|}\}$, the associated vector and diagonal matrix  are $x \in \real^{|\mc I|}$ and $\mathrm{diag}(\{x_i\}_{i\in\mathcal{I}}) \in \real^{|\mathcal{I}|\times|\mathcal{I}|}$.\\
%
Let $\fvec{1}_n$ and $\fvec{0}_n$ be the $n$-dimensional vectors of unit and zero entries. We denote the diagonal vector space $\Span(\fvec 1_{n})$ by $\mathds{1}_{n}$ and its orthogonal complement by $\mathds {1}_n^\perp \triangleq \{x \in \real^n:\, \fvec{1}_n^Tx \!=\! 0\}$. 
%

\emph{Algebraic graph theory: } 
We denote by $G(\mathcal{V},\mathcal{E},A)$ an undirected and weighted graph, where $\mathcal{V}$ is node set, $\mathcal{E} \subseteq \mathcal{V} \times \mathcal{V}$ is the edge set, and $A=A^{T} \in \real^{|\mathcal{V}| \times |\mathcal{V}|}$ is the \emph{adjacency matrix}.
If a number $\ell \in \until{|\mathcal{E}|}$ and an arbitrary direction are assigned to each edge, the \emph{incidence matrix} $B \in \real^{|\mathcal{V}|\times|\mathcal{E}|}$ is defined component-wise as $B_{k\ell} = 1$~if node $k$ is the sink node of edge $\ell$ and as $B_{k\ell} = -1$~if node $k$ is the source node of edge $\ell$; all other elements are zero. 
The \emph{Laplacian matrix} is $L \triangleq B\mathrm{diag}(\{a_{ij}\}_{\{i,j\}\in\mathcal{E}})B^T$.
If the graph is connected, then $\mathrm{ker}(B^T) \!=\! \mathrm{ker}(L) \!=\! \mathds{1}_{|\mc V|}$. For acyclic graphs, $\mathrm{ker}(B) = 
\emptyset$, and for every $x \in \mathds{1}_{|\mathcal{V}|}^\perp$ there is a unique $\xi \in \real^{|\mathcal{E}|}$ satisfying {\em Kirchoff's Current Law} (KCL) $x = B\xi$. 
In a circuit, $x$ are the nodal {current} injections, and $\xi$ are the associated edge flows. 

\emph{Geometry on the $n$-torus}: The set $\mycircle$ denotes the \emph{circle}, an \emph{angle} is a point $\theta \in \mycircle$, and an \emph{arc} is a connected subset of $\mycircle$. Let $|\theta_1 - \theta_2|$ be the \emph{geodesic distance} between two angles $\theta_1,\theta_2 \in \mycircle$. The \emph{$n$-torus} is $\torus^n = \mycircle \times \cdots \times \mycircle$. 
For $\gamma \in [0,\pi/2[$ and a  graph $G(\mathcal{V},\mathcal{E},\cdot)$, let $\overline\Delta_{G}(\gamma) =\{ \theta \in \mathbb T^{|\mathcal{V}|}:\, \max_{\{i,j\} \in \mathcal E} |\theta_{i} - \theta_{j}| \leq \gamma \}$ be the closed set of angle arrays $\theta = ( \theta_{1},\dots,\theta_{n})$ with neighboring angles $\theta_{i}$ and $\theta_{j}$, $\{i,j\} \in \mathcal E$ no further than $\gamma$ apart. Let $\Delta_{G}(\gamma)$ be the interior of $\overline\Delta_{G}(\gamma)$.


\section{Microgrids and their Control Challenges}
\label{Section: Microgrids and their Control Challenges}

\subsection{Microgrids, AC Circuits, and Modeling Assumptions}
\label{Subsection: Review of Microgrids and AC Circuits}

We adopt the standard model of a microgrid as a synchronous linear circuit. The associated connected, undirected, and complex-weighted graph is $G(\mc V,\mc E,A)$ with node set (or buses) $\mc V \!=\! \until n$, edge set (or branches) $\mc E \subset \mc V \!\times\! \mc V$, {impedances $z_{ij} \in \complex$ for each $\{i,j\} \in \mc E$,\,and symmetric weights (or admittances) $a_{ij} = a_{ji} = 1/z_{ij}$.\,\,The admittance matrix $Y \!=\! B \mathrm{diag}( \{ {z^{-1}_{ij}} \}_{\{i,j\}\in\mathcal{E}} )B^T \in \complex^{n \times n}$ is  the associated Laplacian.} 
We restrict ourselves to {\em acyclic}\,(also\,called\,{\em radial}) topologies prevalent in low-voltage distribution networks. 

To each node $i \in \mc V$, we associate an electrical power injection $S_{\mathrm{e},i} = P_{\mathrm{e},i} + \sqrt{-1} Q_{\mathrm{e},i} \in \complex$  and a voltage phasor $V_i = E_{i} e^{\sqrt{-1} \theta_{i}} \in \complex$ corresponding to the magnitude $E_i > 0$ and the phase angle $\theta_i \in \mycircle$ of a harmonic voltage solution to the AC circuit equations.
The complex vector of nodal power injections is then $S_{\mathrm{e}} = V \circ (YV)^C$, where $^{C}$ denotes complex conjugation and $\circ$ is the Hadamard (element-wise) product. 

{
We assume that all lines in the microgrid are made\,\,from the same material, and thus have {\em uniform per-unit-length resistance-to-reactance ratios}. It follows that
$z_{ij} \!=\! |z_{ij}| e^{\sqrt{-1}\varphi}$
for some  angle $\varphi \in {[-\pi/2,\pi/2]}$ for all $\{i,j\} \in \mathcal E$.
The associated admittance matrix is $Y = e^{\sqrt{-1}\left(\pi/2 -\varphi  \right)} \tilde Y$ where
\begin{equation*}
\tilde Y
= { B\,\mathrm{diag}\left( \left\{ - {\sqrt{-1}}/{|z_{ij}|} \right\}_{\{i,j\}\in\mathcal{E}} \right)B^T}\,,
\end{equation*}
is the admittance matrix of the {\em lossless} circuit. After applying the power transformation $\tilde S_{\mathrm{e}} \!=\! S_{\mathrm{e}} e^{\sqrt{-1}(\varphi-\pi/2)}$, the power flow equations are lossless and inductive: $\tilde S_{\mathrm{e}} = V \circ (\tilde Y V)^C$. In the following we assume, without loss of generality, that all lines are purely inductive $\varphi = \pi/2$; see Remark~\ref{Remark: lossy networks} for an additional discussion on this assumption and the power transformation. 

In summary, the active/reactive nodal power injections\,are}%
\begin{subequations}%
\label{eq: power flow}%
\begin{align}%
	P_{\mathrm{e},i} &= \sum\nolimits_{j=1}^{n} \mathfrak{Im}(Y_{ij}) E_{i} E_{j} \sin(\theta_{i} - \theta_{j})
	\,, \;\;\quad i \in \mathcal{V}\,,
	\label{eq: power flow -- active}\\
	\Bigl.
	Q_{\mathrm{e},i} &= - \sum\nolimits_{j=1}^{n} \mathfrak{Im}(Y_{ij}) E_{i} E_{j} \cos(\theta_{i} - \theta_{j})
	\,, \;\; i \in \mathcal{V}.
	\label{eq: power flow -- reactive}
\end{align}
\end{subequations}
We adopt the standard \emph{decoupling approximation} \cite{PK:94,QCZ-TH:13}
 where all voltage magnitudes $E_{i}$ are constant in the active power injections \eqref{eq: power flow -- active} and $P_{\mathrm{e},i} \!=\! P_{\mathrm{e},i}(\theta)$. By continuity and exponential stability, our results are robust to bounded voltage dynamics \cite{JWSP-FD-FB:12u,FD-MC-FB:11v-pnas}, which we illustrate via simulations. 

We partition the set of buses into  loads and inverters, $\mathcal{V} = \mathcal{V}_L \cup \mathcal{V}_I$, and denote their cardinalities by $n \triangleq |\mathcal{V}|$, $n_{L} \triangleq |\mathcal{V}_L|$, and $n_{I} \triangleq |\mathcal{V}_I|$. 
Each load $i \in \mc V_{L}$ demands a constant amount of active power $P_i^* \in \real$ and satisfies the power flow equation
\begin{equation}\label{eq: load flow}
0 = P_i^* - \subscr{P}{e,$i$}(\theta)\,,
\quad i \in \mathcal{V}_L\,.
\end{equation}
{We refer to the buses $\mathcal{V}_L$ strictly as \emph{loads}, with the understanding that they can be either loads demanding constant power $P_i^* <0$, or constant-power sources such as PV inverters performing maximum power point tracking with $P_i^* >0$.}

We denote the rating (maximal power injection) of inverter $i \in \mc V_I$ by $\overline{P}_i \geq 0$. As a necessary feasibility condition, we assume throughout this article that the total load $\sum_{i \in \mathcal{V}_L}P_i^*$ is a net demand {serviceable} by the inverters' maximal generation:
\begin{equation}\label{Eq:LoadRestriction}
	0 \leq - \sum_{i \in \mathcal{V}_L}\nolimits P_i^* \leq \sum\nolimits_{i \in \mathcal{V}_I}\overline{P}_i\,.
\end{equation}
{After appropriate inner control loops are established, an inverter behaves much like a {\em controllable voltage source} behind a reactance \cite{QCZ-TH:13}}, which is the standard model in the literature.
 
\subsection{Primary Droop Control}
\label{Subsection: Review of Droop Control}

The \emph{frequency droop controller} is the main technique for primary control in islanded microgrids \cite{MCC-DMD-RA:93, JMG-JCV-JM-LGDV-MC:11,JAPL-CLM-AGM:06,JMG-MC-TLL-PLC:13,QCZ-TH:13}. 
At inverter $i$, the frequency $\dot\theta_{i}$ is controlled to be proportional to the measured 
 power injection $\subscr{P}{e,$i$}(\theta)$ according to
\begin{equation}\label{eq: primary control}
D_i\dot{\theta}_i = P_i^* - \subscr{P}{e,$i$}(\theta) \,,
\quad i \in \mathcal{V}_I\,,
\end{equation}
where $P_{i}^{*} \in [0,\overline{P}_i]$ is a \emph{nominal injection
  setpoint}, and the proportionality constant $D_{i}\geq 0$ is referred to
as the (inverse) {\em droop coefficient}.
In this notation, $\dot{\theta}_i$ is actually the frequency error $\omega_i - \omega^*$, where $\omega^*$ is the nominal network frequency.

The droop-controlled microgrid is then described by the nonlinear, differential-algebraic equations (DAE)~\eqref{eq: load flow},\eqref{eq: primary control}. 
{
\begin{remark}
\label{Remark: lossy networks}
\textbf{(Droop Controllers for Non-Inductive Networks).}
The droop control equations \eqref{eq: power flow}-\eqref{eq: primary control} are valid for purely inductive lines~without resistive losses. This assumption is typically justified, as the inverter output impedances are controlled to dominate over the network impedances \cite{JMG-LG-JM-MC-JM:05}.\,As discussed prior to equation \eqref{eq: power flow}, this assumption can be made without loss of generality in networks with constant resistance-to-reactance ratios provided that the power injections are transformed as $\tilde S_{\mathrm{e}} = S_{\mathrm{e}} e^{\sqrt{-1}(\varphi-\pi/2)}$, or in components
\begin{equation}
\begin{bmatrix}
\subscr{ \tilde P}{e,$i$} \\  \subscr{\tilde Q}{e,$i$}
\end{bmatrix}
=
\begin{bmatrix}
\sin(\varphi) & - \cos(\varphi) \\ \cos(\varphi) & \sin(\varphi)
\end{bmatrix}
\begin{bmatrix}
 \subscr{P}{e,$i$} \\  \subscr{Q}{e,$i$}
\end{bmatrix}
\,.
\label{eq: power transformation}
\end{equation}
Indeed, \eqref{eq: power transformation} is a common transformation decoupling lossy and lossless injections \cite{JCV-JMG-AL-PR-RT:09,UM-MM:14} 
that is consistent with the fact that active and reactive droop control laws are reversed for resistive lines ($\varphi = 0$) and negated for capacitive lines ($\varphi = -\pi/2$) \cite[Chapter 19.4]{QCZ-TH:13}. Finally, due to continuity and exponential stability, all of our forthcoming results also hold for networks with sufficiently uniform resistance-to-reactance ratios.
\oprocend
\end{remark}
}

%
\subsection{Secondary Frequency Control}
\label{Subsection: Review of Secondary Frequency Control}

The droop controller \eqref{eq: primary control} induces a static error in the steady-state frequency. If the droop-controlled system \eqref{eq: load flow}, \eqref{eq: primary control} settles to a frequency-synchronized solution, $\dot \theta_{i}(t) = \subscr{\omega}{sync} \in \real$ for all $i \in \mc V$, then summing over all equations \eqref{eq: load flow},\eqref{eq: primary control}~yields the synchronous frequency $ \subscr{\omega}{sync}$ as the \emph{scaled power imbalance}
\begin{equation}\label{eq: omega sync}
 \subscr{\omega}{sync}  
 \triangleq \frac{\sum_{i\in\mathcal{V}}P_i^*}{\sum_{i \in \mathcal{V}_I}D_i} \,.
\end{equation}
Notice that $\subscr{\omega}{sync}$ is zero if and only if the nominal injections $P_{i}^{*}$ are balanced: $\sum_{i\in\mathcal{V}}P_i^* = 0$. Since the loads are generally unknown and variable, it is not possible to select the nominal source injections to balance them. Likewise, to render $\subscr{\omega}{sync}$ small, the  coefficients $D_{i}$ cannot be chosen arbitrary large, since the primary control becomes slow and possibly unstable. 

{To eliminate this frequency error, the primary control \eqref{eq: primary control} needs to be augmented with {\em secondary control} inputs $u_{i}(t)$:}
\begin{equation}\label{eq: secondary control}
	D_i\dot{\theta}_i = P_i^* - \subscr{P}{e,$i$}(\theta) + u_{i}(t) \,.
\end{equation}
If there is a synchronized solution to the secondary-controlled equations \eqref{eq: load flow},\eqref{eq: secondary control} with frequency $\subscr{\omega}{sync}^*$ and steady-state secondary control inputs $u_{i}^{*} \!=\! \lim_{t \to \infty} u_{i}(t)$, then we obtain it\,as 
\begin{equation}\label{eq: omega sync with u}
 \subscr{\omega}{sync}^*
 = \frac{\sum_{i\in\mathcal{V}}P_i^*+\sum_{j\in\mathcal{V}_I}u_{i}^{*}}{\sum_{i \in \mathcal{V}_I}D_i} = \omega_{\rm sync} + \frac{\sum_{j\in\mathcal{V}_I}u_{i}^{*}}{\sum_{i \in \mathcal{V}_I}D_i}
 \,.
\end{equation}
Clearly,  there are many choices for the inputs $u_{i}^{*}$ to achieve the control objective $\subscr{\omega}{sync}^* = 0$. However, the inputs $u_i^*$ are typically constrained due to additional performance criteria. 
%

\subsection{Tertiary Operational Control}
\label{Subsection: Review of Tertiary Operational Control}

A tertiary operation and control layer has the objective to minimize an
{\em economic dispatch problem}, that is, an appropriate quadratic cost of the
accumulated~generation: 
\begin{subequations}
\begin{align}
	\minimize_{\theta \in\torus^{n} \,,\, u \in \real^{n_{I}}} 
	& \;\; f(u) = \sum\nolimits_{i \in \mc V_{I}} \frac{1}{2}\alpha_i u_i^2 
	&
	\label{eq: AC optimal dispatch -- cost}\\
	\mbox{subject to} 
	& \qquad P^*_i + u_{i} = \subscr{P}{e,$i$}(\theta) 
	 &  \forall \; i \in \mc V_{I} \,,
	 \label{eq: AC optimal dispatch -- flow generator}\\
	& \qquad P^*_i = \subscr{P}{e,$i$}(\theta)   
	 & \forall \;\, i \in \mc V_{L} \,, 
	 \label{eq: AC optimal dispatch -- flow load}\\
	& \qquad | \theta_{i} - \theta_{j} | \leq \supscr{\gamma}{(AC)}_{ij}
	 & \forall \; \{i,j\} \in \mc E \,, 
	 \label{eq: AC optimal dispatch -- security constraint}\\
	 & \qquad \subscr{P}{e,$i$}(\theta) \in {[0, \overline{P}_i]}
	 & \forall \, i \in \mc V_{I} \,, 
	 \label{eq: AC optimal dispatch -- gen constraint}
\end{align}%
\label{eq: AC optimal dispatch}%
\end{subequations}%
Here, $\alpha_{i} >0$ is the {cost coefficient} for source $i \in \mc V_{I}$.%
\footnote{{See Remark \ref{Remark: Beyond quadratic objective functions} for a discussion of more general objective functions and their implications on the droop curve trading off frequency and active power.}}
The decision variables are the angles $\theta$ and secondary\,control inputs $u$. The non-convex equality constraints \eqref{eq: AC
  optimal dispatch -- flow generator}-\eqref{eq: AC optimal dispatch -- flow load} are~the nonlinear steady-state secondary control equations, the {\em security constraint} \eqref{eq: AC optimal dispatch -- security constraint} limits the power flow on each branch $\{i,j\} \in \mc E$ with $\supscr{\gamma}{(AC)}_{ij} \in {[0,\pi/2[}$, and  \eqref{eq: AC optimal dispatch -- gen constraint} is a~{\em generation constraint}. 
  
{Two typical instances of the economic dispatch \eqref{eq: AC optimal dispatch} are as follows: For $P_{i}^{*} = 0$, $u_{i}$ equals $\subscr{P}{e,$i$}(\theta)$, and the total generation cost is penalized. If the nominal generation setpoints $P_{i}^{*}$ are positive {(for example, scheduled according to a load and renewable generation forecast), 
then $ u_{i}^{*}$ is the operating reserve that must be called upon to meet the real-time net demand.}

\subsection{Heterogeneous Microgrids with Additional Components}
\label{Subsection: Alternative Models of Sources and Loads}

In the following, we briefly list additional components in a microgrid, which can be captured by the model~\eqref{eq: power flow -- active},\eqref{eq: load flow},\eqref{eq: primary control}. 

{\em Synchronous machines:} Synchronous generators (respectively motors) are sources (respectively loads) with dynamics 
\begin{equation}
	M_{i} \ddot \theta_{i} + D_{i} \dot \theta_{i} = P_{i}^{*} - \subscr{P}{e,$i$}(\theta) \,,
	\label{eq: synchronous machine}
\end{equation}
where $M_{i} > 0$ is the inertia, and $D_{i} =  \subscr{D}{diss,$i$} + \subscr{D}{droop,$i$} > 0$ combines dissipation $\subscr{D}{diss,$i$} \dot \theta_{i}$ 
and a droop term $\subscr{D}{droop,$i$} \dot \theta_{i}$ \cite{PK:94}. 
 The constant  injection $P_{i}^{*} \in \real$ is positive for a generator and negative for a load. As shown in \cite[Theorem 5.1]{FD-FB:10w}, the synchronous machine model \eqref{eq: synchronous machine} is {locally} topologically equivalent to a first-order model of the form \eqref{eq: primary control}
{: both models share the same equilibria and the same local  stability\,properties.} 

{\em Inverters with measurement delays:}
The delay between the power measurement $\subscr{P}{e,$i$}(\theta)$ at an inverter $i \in \mc V_{I}$  and the droop control actuation \eqref{eq: primary control} can be explicitly modeled by a first-order lag filter with state $s_{i} \in \real$ and time constant $T_{i}>0$:
\begin{align}\label{eq: primary control with delay}
\begin{split}
		D_i\dot{\theta}_i &= P_i^* - s_{i} \,, \\
		T_{i} \dot s_{i} &= \subscr{P}{e,$i$}(\theta) - s_{i} \,.
\end{split}
\end{align}
As shown in \cite[Lemma 4.1]{JS-DG-JR-TS:13},  after a linear change of variables, the dynamics \eqref{eq: primary control with delay} equal the machine dynamics \eqref{eq: synchronous machine}.

{\em Frequency-dependent loads:} 
If the demand depends on the frequency \cite{NA-SG:13,NA-SG:13b,EM-SHL:14,SY-LC:14,JWSP-FD-FB:12u,CZ-UT-NL-SL:13,NL-CL-ZC-SHL:13}, that is, the left-hand side of \eqref{eq: load flow} is $D_{i} \dot \theta_{i}$ with $D_{i} > 0$, the load dynamics \eqref{eq: load flow} are {formally identical to the} inverter dynamics \eqref{eq: primary control}, {and the microgrid model features no algebraic equations.} This frequency-dependence does not alter the local stability properties \cite{FD-MC-FB:11v-pnas}. 

In summary, all results pertaining to equilibria of the microgrid model~\eqref{eq: load flow},\eqref{eq: primary control} and their local stability extend to synchronous machines, inverters with measurement delays, and  frequency-dependent loads. Likewise, all secondary or tertiary control strategies can be equally applied. With these extensions in mind, we focus on the  microgrid model~\eqref{eq: load flow},\eqref{eq: primary control}.
%
%
%
\section{Decentralized Primary Control Strategies}
\label{Section: Primary Control Strategies}

In this section, we study the fundamental properties of the droop-controlled microgrid \eqref{eq: load flow},\eqref{eq: primary control}. In Section \ref{Section: Secondary Control Strategies}, we design appropriate secondary controllers, which preserve the properties of primary control even if the load profile is unknown.
%
\subsection{Symmetries, Synchronization, and Transformations}
\label{Subsection: Symmetries, Synchronization, and Transformations}

{We begin our analysis by reviewing the symmetries of the droop-controlled microgrid \eqref{eq: load flow},\eqref{eq: primary control}.
Observe that if the system \eqref{eq: load flow},\eqref{eq: primary control} possesses a stable and synchronized solution with frequency $\subscr{\omega}{sync}$ as given in \eqref{eq: omega sync}, then it possesses stable equilibria in a rotating coordinate frame with frequency $\subscr{\omega}{sync}$. The effect of secondary control \eqref{eq: secondary control} regulating the frequency can be understood as carrying out such a coordinate transformation to an appropriately rotating frame. In this subsection, we formally establish the equivalence of these three ideas so that\,we\,can restrict our attention to a shifted system in rotating\,coordinates.}

The microgrid equations \eqref{eq: load flow},\eqref{eq: primary control} feature an inherent {\em rotational symmetry}: they are invariant under a rigid rotation of all angles. Formally, let $\rot_s(r)\in\mycircle$ be the rotation of a point $r \in \mycircle$
counterclockwise by the angle $s\in{[0,2\pi]}$.  For $(r_{1},\dots,r_{n}) \in \torus^n
$, define the equivalence class
\begin{equation*}
  [(r_1,\dots,r_n)] \!=\! \left\{ ( \rot_s(r_1), \dots, \rot_s(r_n) ) \!\in\! \torus^{n}:\, \! s \in {[0,2\pi]} \right\}.
\end{equation*}
Thus, a synchronized solution $\theta^{*}(t)$ of \eqref{eq: load flow},\eqref{eq: primary control} is part of a one-dimensional connected {\em synchronization manifold} $[\theta^{*}]$. 
For $\subscr{\omega}{sync} = 0$, a synchronization manifold is also an equilibrium manifold of \eqref{eq: load flow},\eqref{eq: primary control}.
 In the following, when we refer to a synchronized solution as ``stable'' or ``unique'', these properties are to be understood 
 modulo rotational symmetry.

Recall that, without secondary control, the synchronous~frequency $ \subscr{\omega}{sync}$ is the {scaled power imbalance} \eqref{eq: omega sync}. By transforming to a rotating coordinate frame 
with frequency $\subscr{\omega}{sync}$, {that is, $\theta_{i}(t) \mapsto \rot_{\subscr{\omega}{sync} t}(\theta_{i}(t))$ (with slight abuse of notation, we maintain the variable $\theta$),} 
a synchronized solution of \eqref{eq: load flow},\eqref{eq: primary control} is equivalent to an equilibrium of the \emph{shifted control system}
\begin{subequations}
\begin{align}
	0 &= \widetilde P_{i} - P_{\mathrm{e},i}(\theta)
	\,, &\quad i \in \mc V_{L} \,,
	\label{eq: primary-controlled system - transf - load}\\
	D_{i} \dot \theta_{i} &= \widetilde P_{i} - P_{\mathrm{e},i}(\theta)
	\,, &\quad i \in \mc V_{I} \,,
	\label{eq: primary-controlled system - transf - inv}
\end{align}%
\label{eq: primary-controlled system}%
\end{subequations}%
where the {\em shifted power injections} are $\widetilde P_{i} = P_{i}^{*}$ for $i \in \mc V_{L}$, and $\widetilde P_{i} = P_{i}^{*} - D_{i}\subscr{\omega}{sync}$ for $i \in \mc V_{I}$. We emphasize that the shifted injections in \eqref{eq: primary-controlled system} are balanced: $\widetilde P \in \mathds 1_{n}^\perp$.  
%
Notice that, equivalently to transforming to a rotating frame with frequency $\subscr{\omega}{sync}$ (or replacing $P$ by $\tilde P$), we can assume that the secondary control input in \eqref{eq: load flow},\eqref{eq: secondary control}  takes the constant value $u_{i} = -D_{i}\subscr{\omega}{sync}$ for all $i \in \mc V_{I}$ to arrive at the shifted control system \eqref{eq: primary-controlled system}. 

We summarize these observations in the following lemma.
\begin{lemma}\textbf{(Synchronization Equivalences).}
\label{Lemma: Synchronization Equivalences}
The following statements are equivalent:
\begin{enumerate}

	\item The primary droop-controlled microgrid \eqref{eq: load flow},\eqref{eq: primary control} possesses a locally exponentially stable and unique 
synchronization manifold $t \mapsto [\theta(t)] \subset \torus^{n}$ for all $t \geq 0$;
	
    \item The secondary droop-controlled microgrid \eqref{eq: load flow},\eqref{eq: secondary control} with constant secondary-control input $u_{i} = -D_{i}\subscr{\omega}{sync}$ for all $i \in \mc V_{I}$ possesses a locally exponentially stable and unique 
equilibrium manifold $[\bar \theta] \subset \torus^{n}$;

	\item The shifted control system \eqref{eq: primary-controlled system} possesses a locally exponentially stable and unique 
equilibrium $[\tilde \theta] \subset \torus^{n}$.

\end{enumerate}
If the equivalent statements (i)-(iii) are true, then all systems have the same synchronization manifolds $[\theta(t)] = [\bar\theta] = [\tilde\theta]\subset \torus^{n}$ and the same power injections $\subscr{P}{e}(\theta(t)) = \subscr{P}{e}(\bar\theta) = \subscr{P}{e}(\tilde\theta)$. Additionally, $\theta(t) = \rot_{\subscr{\omega}{sync}t}(\xi_{0})$ for some $\xi_{0} \in [\bar\theta] = [\tilde\theta]$.
\end{lemma}

In light of Lemma \ref{Lemma: Synchronization Equivalences}, we restrict the {forthcoming} discussion in this section to the shifted control system \eqref{eq: primary-controlled system}.

Observe also that equilibria of \eqref{eq: primary-controlled system} are invariant under {\em constant scaling} of all droop coefficients: if $D_{i}$ is replaced by $D_{i} \cdot \beta$ for some $\beta \in \real$, then ${\subscr{\omega}{sync}}$ changes to $\subscr{\omega}{sync}/\beta$. Since the product $D_{i} \cdot \subscr{\omega}{sync}$ remains constant, the equilibria of \eqref{eq: primary-controlled system} do not change. Moreover, if $\beta>0$, then the stability properties of  equilibria do not change since time can be rescaled as $t \mapsto t / \beta$.

\subsection{Existence, Uniqueness, \& Stability of Synchronization}
\label{Subsection: Basic Stability Theorem}

In vector form, the equilibria of \eqref{eq: primary-controlled system} satisfy
\begin{equation}
	\widetilde P = B \mc A \sinbf(B^{T} \theta)
	\,,
	\label{eq: equilibria - primary control - 1}
\end{equation}
where $\mc A = \diag(\{ \mathfrak{Im}(Y_{ij}) E_{i} E_{j} \}_{\{i,j\} \in \mc E})$ and $B \in \real^{|\mc V| \times |\mc E|}$ is the incidence matrix. 
For an acyclic network, $\mathrm{ker}(B) = \emptyset$, and the unique vector of branch flows $\xi \in \real^{|\mathcal{E}|}$ (associated to the shifted injections $\widetilde{P}$) is given by KCL as $\xi =  B^\dagger \widetilde{P} = (B^TB)^{-1}B^T \widetilde{P}$. 
Hence, equations \eqref{eq: equilibria - primary control - 1} equivalently read~as
\begin{equation}
	\xi = \mc A \sinbf(B^{T} \theta)
	\label{eq: equilibria - primary control - 2}
	\,.
\end{equation}	
Due to boundedness of the sinusoid, a necessary condition for solvability of equation \eqref{eq: equilibria - primary control - 2} is  $\| \mc A^{-1} \xi \|_{\infty}<1$. The following result shows that this condition is also sufficient and guarantees stability of an equilibrium manifold of \eqref{eq: primary-controlled system} \cite[Theorem 2]{JWSP-FD-FB:12u}.

\begin{theorem}\textbf{(Existence and Stability of Synchronization).} 
\label{Thm:Stab}
Consider the shifted control system \eqref{eq: primary-controlled system}. Let $\xi \in \real^{|\mathcal{E}|}$ be the unique vector of power flows satisfying the KCL, given by $\xi = B^\dagger \widetilde{P}$\@. The following two statements are equivalent:
\begin{enumerate}

\item[(i)] \textbf{\em Synchronization: } there exists an arc length $\gamma \in [0,\pi/2[$ such that the shifted control system \eqref{eq: primary-controlled system} possesses a locally exponentially stable and unique 
equilibrium manifold $[\theta^*] \subset \overline\Delta_{G}(\gamma)$;

\item[(ii)] \textbf{\em Flow feasibility:} the power flow is feasible, that is, 
\begin{equation} \label{eq: sync condition}
\Gamma \triangleq \|\mathcal{A}^{-1}\xi\|_{\infty} < 1.
\end{equation}
\end{enumerate}
%
If the equivalent statements (i) and (ii) hold true, then the quantities $\Gamma \in [0,1[$ and $\gamma \in [0,\pi/2[$ are related uniquely via $\Gamma = \sin(\gamma)$, and $\boldsymbol{\sin}(B^T\theta^*) = \mc A^{-1}\xi$\@.
\end{theorem}

\subsection{Power Flow Constraints and Proportional Power Sharing}
\label{Subsection: Proportional Power Sharing}

{While Theorem \ref{Thm:Stab} gives the exact stability condition, it offers no guidance on how
to select the control parameters $(P_i^*,D_i)$ to
achieve a set of desired steady-state power injections.
One desired objective is that all sources meet their the actuation
constraints $\subscr{P}{e,$i$}(\theta) \in [0,\overline{P}_i]$ and share the {load} in a fair way according to their power ratings \cite{QCZ-TH:13, JMG-JCV-JM-LGDV-MC:11,JMG-MC-TLL-PLC:13}.}

\begin{definition}\textbf{{(Proportional Power Sharing).}}\label{Def:PropShar}
Consider an equilibrium manifold $[\theta^{*}] \subset \torus^{n}$ of the shifted control system \eqref{eq: primary-controlled system}. The inverters $\mc V_{I}$ share the total load $\sum_{i \in \mathcal{V}_L}P_i^*$ {\em proportionally according to their power ratings} if for all $i,j \in \mathcal{V}_I$
\begin{equation}
	{\subscr{P}{e,$i$}(\theta^{*})}/{\overline{P}_i} = {\subscr{P}{e,$j$}(\theta^{*})}/{\overline{P}_j} 
	\,.
	\label{eq: prop power sharing}
\end{equation}
\end{definition}

We also define a useful choice of droop coefficients.
\begin{definition}\textbf{{(Proportional Droop Coefficients and Nominal Injection Setpoints).}}\label{Def:Propor}
The droop coefficients and injections setpoints are selected proportionally if for all $i,j\!\in\! \mathcal{V}_I$
\begin{equation}
	P_i^*/D_i = P_j^*/D_j 
	\,\mbox{ and }\,
	P_i^*/\overline{P}_i = P_j^*/\overline{P}_j
	\,.
	\label{eq: prop choice of droop coeff}
\end{equation}
\end{definition}

A proportional choice of droop control coefficients leads\,to\,a fair load  sharing among the inverters according to their ratings and subject to their actuation constraints -- a result that also holds for lossy and meshed circuits \cite[Theorem 7]{JWSP-FD-FB:12u}:
\begin{theorem}\label{Thm:PowerFlowConstraintsTree}
\textbf{{(Power Flow Constraints and Power Sharing).}}
Consider an equilibrium manifold $[\theta^{*}] \subset \torus^{n}$ of the shifted control system \eqref{eq: primary-controlled system}. Let the droop coefficients be selected proportionally. The following statements are equivalent:

\begin{enumerate}
\item[(i)] {\textbf{\em Injection constraints: }} $0 \leq \subscr{P}{e,$i$}(\theta^{*}) \leq \overline{P}_i$, $\,\,\,\forall i \in \mathcal{V}_I$;
\item[(ii)] {\textbf{\em Serviceable load: }} $0 \leq - \sum_{i \in \mathcal{V}_L}P_i^* \leq \sum_{j \in \mathcal{V}_I}\overline{P}_j\,.$
\end{enumerate}
Moreover, the inverters share the total load $\sum_{i \in \mathcal{V}_L}P_i^*$ proportionally according to their power ratings. 
\end{theorem}

%
\subsection{Power Flow Shaping}
\label{Subsection: Flow Shaping}

We now address the following ``{reachability}'' question: given a set of desired power injections for the inverters, can one select the droop coefficients to generate these injections?

We define a {\em power injection setpoint} as a point of power balance,
at fixed load demands and subject to the basic feasibility condition \eqref{eq: sync condition} given in
Theorem \ref{Thm:Stab}.

\begin{definition}\textbf{(Feasible Power Injection Setpoint).}\label{Def:Inj}
Let $\gamma \in [0,\pi/2[$. A vector $P^{\rm set} \in \real^n$ is a \emph{$\gamma$-feasible power injection setpoint} if it satisfies the following three properties:
\begin{enumerate}

	\item \textbf{\em Power balance:} $P^{\rm set} \in \mathds{1}_n^{\perp}$;
	
	\item \textbf{\em Load invariance:}  $\supscr{P}{set}_{i} = P_i^*$ for all loads $i \in \mathcal{V}_L$;
	
	\item \textbf{\em $\gamma$-feasibility:} the associated branch  power flows $\xi^{\rm set} = B^\dagger P^{\rm set}$ are feasible, that is, $\|\mathcal{A}^{-1}\xi^{\rm set}\|_\infty \leq \sin(\gamma)$.
	

\end{enumerate}
\end{definition}

\smallskip

The next result characterizes the relationship between droop controller
designs and $\gamma$-feasible injection setpoints. For simplicity, we omit
the singular case where $\omega_{\rm sync} = 0$, since in this case the
droop coefficients offer no control over the steady-state inverter 
injections $P_{\mathrm{e},i}(\theta^{*}) = P_i^* - D_i\omega_{\rm sync}$.

\begin{theorem}
\textbf{(Power Flow Shaping).}
\label{Theorem: power injection setpoint design}
Consider the shifted control system \eqref{eq: primary-controlled system}. Assume $\omega_{\rm sync} \neq 0$, let $P^{\rm set} \in \mathds 1_{n}^{\perp}$, and let $\gamma \in [0,\pi/2[$. The following statements are equivalent:
\begin{enumerate}
\item[(i)] \textbf{\em Coefficient selection: }there exists a selection of droop coefficients $D_{i}$, $i \in \mathcal{V}_I$, such that the steady-state injections satisfy $P_{\rm e}(\theta^*) = P^{\rm set}$, with $[\theta^*]\subset \overline\Delta_{G}(\gamma)$; 
\item[(ii)] \textbf{\em Feasibility: }$P^{\rm set}$ is a $\gamma$-feasible power injection setpoint.
\end{enumerate}
If the equivalent statements (i) and (ii) hold true, then the quantities $D_{i}$ and $\supscr{P}{set}_{i}$ are related with arbitrary $\beta \neq 0$ as
\begin{equation}\label{eq: DroopSetpoint}
D_i = \beta ({P}_i^{*} - \supscr{P}{set}_{i}) \,,\quad i \in \mathcal{V}_I \,.
\end{equation}
Moreover, $[\theta^*]$ is locally exponentially stable if and only if $\beta({P}_i^{*} - \supscr{P}{set}_{i})$ is nonnegative for all $i \in \mathcal{V}_I$. 
\end{theorem}

\begin{IEEEproof}
\textbf{(i)}$\implies$\textbf{(ii): }
Since $\theta^* \in \overline\Delta_{G}(\gamma)$ and $P_{\rm e}(\theta^{*}) \in \mathds{1}_n^\perp$, Theorem \ref{Thm:Stab} shows that $P^{\rm set}$ is a $\gamma$-feasible injection setpoint. 

\textbf{(ii)}$\implies$\textbf{(i): }Let $P^{\rm set}$ be a $\gamma$-feasible injection setpoint. 
Consider the droop coefficients $D_i = \beta (P_i^* - \supscr{P}{set}_{i})$. Since $\omega_{\rm sync} \neq 0$, for each $i \in \mathcal{V}_I$ we obtain the steady-state injection
\begin{align*}
P_{\mathrm{e},i}(\theta^{*}) &= \widetilde{P}_i = P_i^* - D_i\subscr{\omega}{sync}\\
%
%
%
&= P_i^* - \beta (P_i^*-\supscr{P}{set}_{i})\frac{1}{\beta}\underbrace{\frac{\sum_{i \in \mathcal{V}}P_i^*}{\sum_{i \in \mathcal{V}_I} (P_i^*-P_i^{\rm set})}}_{=1}%
%
%
%
%
%
%
= \supscr{P}{set}_{i}\,,
\end{align*}
where we used $\sum_{i \in \mathcal{V}_I}P_i^{\rm set} = -\sum_{i \in \mathcal{V}_L}P_{i}^{*}$. Since $P_{\mathrm{e},i}(\theta^{*}) \!=\! P_i^* \!=\! P_{i}^{\rm set}$ for each $i \in \mathcal{V}_L$, we have $P_{\rm e}(\theta) = P^{\rm set}$. Since $P^{\rm set}$ is $\gamma$-feasible, $\theta^{*}$ is well defined in $ \overline\Delta_{G}(\gamma)$.
By the reasoning leading to Theorem \ref{Thm:Stab} (see \cite[Theorem 2]{JWSP-FD-FB:12u}), the shifted system \eqref{eq: primary-controlled system} is stable if and only if all $D_i$ are nonnegative. 
\end{IEEEproof}

\begin{remark}\textbf{(Generation Constraints).}
For a $\gamma$-feasible injection setpoint, the  inverter
generation constraint $P_{i}^{\rm set} \in [0,\overline{P}_i]$ is generally not {guaranteed to be met}. This constraint is feasible if $P_i^* = 0$ ($i \in
\mathcal{V}_I$) and an additional  condition\,holds:
\begin{equation}\label{Eq:JohnBound}
-\sum_{j \in \mathcal{V}_L}\nolimits P_j^* \leq \left({\overline{P}_i}/{D_i}\right) \sum_{j \in \mathcal{V}_I}\nolimits D_j\,,\quad i \in \mathcal{V}_I\,.
\end{equation}
The inequalities \eqref{Eq:JohnBound} limit the heterogeneity of the inverter power injections, and are sufficient for the load serviceability condition \eqref{Eq:LoadRestriction}, as one can see by rearranging and summing over all $i \in \mathcal{V}_I$. A similar result holds for the choice $P_i^* = \overline{P}_i$.
\oprocend
\end{remark}


\section{Centralized, Decentralized, and Distributed Secondary Control Strategies}
\label{Section: Secondary Control Strategies}
%

The primary droop controller \eqref{eq: primary control} results in the static frequency error $\subscr{\omega}{sync}$ in \eqref{eq: omega sync}. The purpose of the secondary control $u_i(t)$ in \eqref{eq: secondary control} is to eliminate this frequency error despite unknown and variable loads. In this section, we investigate different decentralized and distributed secondary control strategies. 

\subsection{Decentralized Secondary Integral Control}
\label{Subsection: Decentralized Secondary Integral Control}

To investigate decentralized secondary control, we partition the set of inverters as $\mc V_{I} = \mc V_{I_{P}} \cup \mc V_{I_{S}}$, where the action of the $\mathcal{V} _{I_{P}}$ inverters is restricted to primary droop control, and the $\mathcal{V} _{I_{S}}$ inverters use the local frequency error for integral control:
\begin{equation}
\begin{aligned}
		u_{i}(t) &=  - p_{i} \;\;,\;\; k_{i} \dot p_{i} =  \dot \theta_i
		\,,\qquad i \in \mc V_{I_{S}} \,,\\
		u_i(t) &= 0\,,\qquad \qquad \qquad \qquad \;\, i \in \mc V_{I_P}\,.
\end{aligned}
	\label{eq: decentralized control}
\end{equation}
{Consider the case $|\mc V_{I_{S}}|=1$, which mimics AGC inside a control area of a transmission network. It can be shown, as a direct corollary to Theorem \ref{Theorem: Stability of Partial Secondary Control} (in Section \ref{Subsection: Partial Secondary Control}), that this controller achieves frequency regulation but fails to maintain the power sharing. Additionally, if a steady-state exists, $p_i$ must converge to the total power imbalance $\sum_{i \in \mathcal{V}}P_i^*$, which places a large and unpredictable burden on a single generator.}

{For $|\mc V_{I_{S}}|>1$, it is well-known in control \cite{KJA-TH:06} and in power systems \cite{PK:94}, that multiple integrators in an interconnected system lead undesirable properties such as undesirable equilibria.
\begin{lemma}[\bf Steady-states of decentralized integral control]
\label{Lemma: Steady-states of decentralized integral control}
Consider the droop-controlled microgrid \eqref{eq: load flow},\eqref{eq: secondary control} with decentralized secondary integral control \eqref{eq: decentralized control} at nodes $\mathcal{V}_{I_{S}}$.  If the system reaches an equilibrium, the steady-state source injections are determined only up to a subspace of the same dimension as the number of decentralized integrators $|\mathcal{V}_{I_{S}}|$.
\end{lemma}

\begin{IEEEproof}
We partition the injections according to loads $\mc V_{L}$, inverters $\mc V_{I_{P}} = \mc V_{I}\setminus \mathcal{V}_{I_S}$ with only primary control \eqref{eq: primary control}, and inverters $\mathcal{V}_{I_S}$ with decentralized integral control \eqref{eq: secondary control},\eqref{eq: decentralized control}
 as $\widetilde P =
(\widetilde P_{L},\widetilde P_{I_{P}},\widetilde P_{I_{S}})$ and 
$P_{\mathrm{e}}(\theta) = (P_{\mathrm{e},L}(\theta),
P_{\mathrm{e},I_{P}}(\theta),P_{\mathrm{e},I_{S}}(\theta))$ and obtain the closed-loop equilibria from the equations
  \begin{equation*}
  \begin{bmatrix}
	\boldsymbol{0} \\ \boldsymbol{0} \\ \boldsymbol{0} \\ \boldsymbol{0}
	\end{bmatrix}
	=
		\begin{bmatrix}
		I & \boldsymbol{0} & \boldsymbol{0} & \boldsymbol{0} \\ 
		\boldsymbol{0}  &  I & \boldsymbol{0}  & \boldsymbol{0}  \\
		\boldsymbol{0} & \boldsymbol{0}  &  I & \fvec 1 \\
		\boldsymbol{0} &  \boldsymbol{0}   & \boldsymbol{0}  & \boldsymbol{0} 
	\end{bmatrix}
	\begin{bmatrix}
		\widetilde P_{L} -  P_{\mathrm{e},L}(\theta)  \\ \widetilde P_{I_{P}} -  P_{\mathrm{e},I_{P}}(\theta) \\ \widetilde P_{I_{S}} -  P_{\mathrm{e},I_{S}}(\theta) \\ -p
	\end{bmatrix}
	.
  \end{equation*}

  It follows that the matrix determining the inverter equilibrium injections $P_{\mathrm{e},I_{S}}(\theta)$ has a $|\mathcal{V}_{I_{S}}|$-dimensional nullspace.
\end{IEEEproof}

Lemma~\ref{Lemma: Steady-states of decentralized integral control} implies that the primary objectives (such as proportional power sharing or power flow shaping) cannot be recovered and the steady-state injections depend on initial values, loads, and exogenous disturbances. These steady-state}
%
 subspaces  correspond to~different choices of $u^{*}$ rendering $\subscr{\omega}{sync}^*$ to zero in \eqref{eq: omega sync with u}. 
One way to remove these undesirable subspaces is to implement \eqref{eq: decentralized control} via the low-pass filter
\begin{equation}
		u_{i}(t) =  - p_{i} \;\;,\;\; k_{i} \dot p_{i} =  \dot \theta_i - \epsilon p_{i}
		\,,\qquad i \in \mc V_{I_{S}} \,,
	\label{eq: leaky centralized control}
\end{equation}
For small $\epsilon>0$ and large $k>0$ (enforcing a time-scale separation), the controller \eqref{eq: leaky centralized control} achieves practical stabilization but does not exactly regulate the frequency \cite{NA-SG:13b}.
In conclusion, the decentralized control \eqref{eq: decentralized control} and its variations generally fail to achieve fast frequency regulation while maintaining power sharing among generating units.
Additionally, a single microgrid source may not have the authority or the capacity to perform secondary control.
In the following, we analyze distributed strategies that exactly recover the primary injections.
%
\subsection{Centralized Averaging PI (CAPI) Control}
\label{Subsection: Centralized Averaging PI Control}

Different distributed secondary control strategies have been proposed in \cite{QS-JG-JMV:13-updated,HL-BJC-WZ-XS:13}. In \cite{QS-JG-JMV:13-updated}, an integral feedback of a {\em weighted average frequency} among all inverters is proposed:%
\footnote[3]{{Aside from the integral feedback \eqref{eq: CAPI},} the controller in \cite{QS-JG-JMV:13-updated} also contains a proportional feedback of the average frequency. We found that such a proportional feedback does generally not preserve the equilibrium injections, and we omit it here. The controller in \cite{QS-JG-JMV:13-updated} uses an arithmetic average with all $D_{i}=1$ in \eqref{eq: CAPI}. Since the synchronization frequency \eqref{eq: omega sync} is obtained by a weighted average and since $D_{i}\dot\theta_{i}$ is the inverter injection $P_{\mathrm{e},i}(\theta)$ (for $P_{i}^{*}=0$),\,we found the choice \eqref{eq: CAPI} more natural. Simulations indicate that any convex combination of the inverter frequencies yields identical results.}
\begin{equation}
		u_{i}(t) =  - p_{i} \;\;,\;\; k_{i} \dot p_{i} =  \frac{\sum\nolimits_{j \in \mc V_{I}} D_{j}\dot \theta_{j}}{{\sum_{j \in \mathcal{V}_I}D_j}}	
	\,, \quad i \in \mc V_{I} \,,
	\label{eq: CAPI}
\end{equation}
Here, $p_i$ is the secondary variable and $k_i > 0$. For  $P_{i}^{*} = 0$, the average frequency in \eqref{eq: CAPI} is the sum of the inverter injections $P_{\mathrm{e},i}(\theta)$. In this case, \eqref{eq: CAPI} equals the secondary control strategy in \cite{HL-BJC-WZ-XS:13}, where the averaged inverter injections are integrated.

By counter-examples, it can be shown that the 
secondary controller \eqref{eq: CAPI} does not have the power sharing property of the shifted control system \eqref{eq: primary-controlled system} unless the values of $D_{i}$ and $k_{i}$ are carefully tuned.
In the following, we suggest the choice
\begin{equation}
	k_{i} = k / D_{i} \,, \quad i \in \mc V_{I} 
	\label{eq: CAPI - tuning}
	\,,
\end{equation}
where $k > 0$. 
The closed loop \eqref{eq: load flow}, \eqref{eq: secondary control}, \eqref{eq: CAPI}, \eqref{eq: CAPI - tuning} is given by%
 \begin{subequations}%
 \label{eq: secondary control -- QS}%
 \begin{align}
 		0 &= P_i^* - \subscr{P}{e,$i$}(\theta)
	\,,\quad & i \in \mathcal{V}_L  \,,
	\label{eq: secondary control -- QS - 0} \\
	D_{i} \dot \theta_{i} &= P_{i}^{*} - \subscr{P}{e,$i$}(\theta) - p_i
	\,, \label{eq: secondary control -- QS - 1} \quad & i \in \mc V_{I} \,,\\
	k \dfrac{\dot {p}_{i}}{D_{i}} &= \frac{\sum\nolimits_{j \in \mc V_{I}} D_{j}\dot \theta_{j}}{{\sum_{j \in \mathcal{V}_I}D_j}}
	\,, \quad & i \in \mc V_{I} \,.
	\label{eq: secondary control -- QS - 2}
 \end{align}
 \end{subequations}
By changing variables $q_{i} \triangleq p_{i}/D_{i} - \subscr{\omega}{sync}$ for $i \in \mc V_{I}$ and observing that $k \dot {q}_{i} = {\sum\nolimits_{j \in \mc V_{I}} D_{j}\dot \theta_{j}}/{{\sum_{j \in \mathcal{V}_I}D_j}}$ is identical for all $i \in  \mc V_{I} $, we can rewrite the closed-loop equations \eqref{eq: secondary control -- QS} as
  \begin{subequations}
   \begin{align}
   		0 &= \widetilde P_i - \subscr{P}{e,$i$}(\theta)
	\,,\quad & i \in \mathcal{V}_L  \,,
	\label{eq: secondary control -- QS_transf - 0} \\
	D_{i} \dot \theta_{i} &= \widetilde P_{i} - \subscr{P}{e,$i$}(\theta)  - D_{i} q 
	\,,\label{eq: secondary control -- QS_transf - 1}  \quad & i \in \mc V_{I} \,, \\
	k \dot {q} &= \frac{\sum\nolimits_{j \in \mc V_{I}} D_{j}\dot \theta_{j}}{{\sum_{j \in \mathcal{V}_I}D_j}}
	 \,,
	\label{eq: secondary control -- QS_transf - 2}
   \end{align}%
     \label{eq: secondary control -- QS_transf}%
 \end{subequations}%
 where $\widetilde P_{i}$ is as in \eqref{eq: primary-controlled system}.
In this transformed system, \eqref{eq: secondary control -- QS_transf - 2} can be implemented as a {\em centralized controller}: {it receives frequency measurements from all inverters and broadcasts the secondary control variable $q$ back to all units.} Due to this insight on the communication complexity, we refer to \eqref{eq: CAPI}-\eqref{eq: CAPI - tuning} as {\em
   centralized averaging proportional integral} (CAPI) control.

 \begin{theorem}\textbf{(Stability of CAPI-Controlled Network).}
\label{Theorem: Stability of CAPI}
Consider the droop-controlled microgrid \eqref{eq: load flow},\eqref{eq: secondary control} with $P_{i}^{*} \in [0,\overline{P}_i]$ and $D_{i}>0$ for $i \in \mathcal{V}_I$. Assume a complete communication topology among the inverters $\mc V_{I}$, and let $u_i(t)$ be given by the CAPI controller \eqref{eq: CAPI} with the parametric choice \eqref{eq: CAPI - tuning}. The following two statements are equivalent:
\begin{enumerate}
\item[(i)] \textbf{\em Stability of primary droop control:  } the droop control stability condition \eqref{eq: sync condition} holds;

\item[(ii)] \textbf{\em Stability of CAPI control: } there exists an arc length $\gamma \in [0,\pi/2[$ such that the closed loop \eqref{eq: secondary control -- QS} possesses a locally exponentially stable and unique 
equilibrium manifold $([\theta^*],p^{*}) \subset \overline\Delta_{G}(\gamma)\times \real^{n_{I}}$.
\end{enumerate}
If the equivalent statements (i) and (ii) are true, then $[\theta^{*}]$ {and all injections} are as in Theorem \ref{Thm:Stab}, and $p_{i}^{*} = D_{i} \subscr{\omega}{sync}$, $i \in \mathcal{V}_{I}$.
\end{theorem}

 \begin{IEEEproof}
{We start by writing system \eqref{eq: secondary control -- QS_transf} in
vector form. Let $D_{I} = \diag(\{D_{i}\}_{i \in \mc V_{I}})$, let
$\subscr{D}{tot} \!=\! \sum_{i \in \mathcal{V}_I}D_i$,
and partition the angles according to loads $\mc V_{L}$ and inverters $\mc V_{I}$ as $\theta =
(\theta_{L},\theta_{I})$.  With this notation, the closed loop \eqref{eq:
  secondary control -- QS_transf} reads in vector form\,as}
  \begin{multline}
   	\underbrace{
 	\begin{bmatrix}
	I & \boldsymbol{0} & \boldsymbol{0} \\ \boldsymbol{0} & \!\!D_{I}\!\! & \boldsymbol{0} \\ \boldsymbol{0} & \boldsymbol{0} & \!\!k \cdot \subscr{D}{tot}
	\end{bmatrix}
	}_{\triangleq Q_{1}}
	\begin{bmatrix}
	\boldsymbol{0} \\ \dot \theta_{I} \\ \dot q
	\end{bmatrix}
	=
		\begin{bmatrix}
		I & \boldsymbol{0} & \boldsymbol{0} \\ 
		\boldsymbol{0} & I & \!D_{I} \fvec 1\! \\
		\boldsymbol{0} &  \!\!\fvec 1^{T} \!\!  & \!\subscr{D}{tot}\!
	\end{bmatrix}
	\begin{bmatrix}
		\!\widetilde P_{L} -  P_{\mathrm{e},L}(\theta)\!  \\ \!\widetilde P_{I} -  P_{\mathrm{e},I}(\theta)\! \\ -q
	\end{bmatrix}
	\\
	=
 		\underbrace{
	\begin{bmatrix}
	I & \boldsymbol{0} & \boldsymbol{0}\\
	\boldsymbol{0} & \!\!D_{I}\!\! & \boldsymbol{0} \\
	\boldsymbol{0} & \boldsymbol{0} & 1
	\end{bmatrix}
	}_{\triangleq Q_{2}}
	\underbrace{
	\begin{bmatrix}
		I & \boldsymbol{0} & \boldsymbol{0} \\ 
		\boldsymbol{0} & D_{I}^{-1} & \fvec 1 \\
		\boldsymbol{0} &  \!\!\fvec 1^{T}\!\!  & \!\subscr{D}{tot}
	\end{bmatrix}
	}_{\triangleq Q_{3}}
	\underbrace{
	\begin{bmatrix}
		\!\widetilde P_{L} -  P_{e,L}(\theta)\!  \\ \!\widetilde P_{I} -  P_{e,I}(\theta)\! \\ -q
	\end{bmatrix}
	}_{\triangleq x}
	\,.
	 \label{eq: closed-loop dynamics -- QS}
 \end{multline}
 The matrices $Q_{1}$ and $Q_{2}$ are nonsingular, while $Q_{3}$ is singular with $\textup{ker}(Q_{3}) = \Span([ 0 \,,\,D_{I} \fvec 1_{n_{I}} \,,\, -1])$.
 On the other hand, 
$
[ \fvec 1_{n}^T \; \fvec 1_{n}^T \; 0 ] x = 0
 $
due to balanced injections $\fvec 1^{T} \widetilde P = 0$ and flows $\fvec 1^{T} P_{\mathrm{e}}(\theta) = 0$. It follows that $x \not\in \textup{ker}(Q_{3})$. Thus, possible equilibria of \eqref{eq: closed-loop dynamics -- QS} are $x = \fvec 0_{n}$, that is, the equilibria $[\theta^{*}]$ from \eqref{eq: primary-controlled system} and $q^{*} = 0$. By Theorem \ref{Thm:Stab}, the equation $x= \fvec 0_{n}$ is solvable for $[\theta^*] \in \overline\Delta_{G}(\gamma)$ if and only if condition \eqref{eq: sync condition} holds.

To establish stability, observe that the negative power flow Jacobian $-\partial P_{\mathrm{e}}(\theta) /\partial \theta $ equals the Laplacian matrix ${\mc L}(\theta) \!=\! B \diag(\{a_{ij}\}_{\{i,j\} \in \mc E})B^{T}$ with $a_{ij} \triangleq \mathfrak{Im}(Y_{ij}) E_{i} E_{j} \cos(\theta_{i}-\theta_{j})$ as weights  \cite[Lemma 2]{FD-MC-FB:11v-pnas}.
For $[\theta^{*}] \in \overline\Delta_{G}(\gamma)$, all weights $a_{ij}>0$ are strictly positive for $\{i,j\} \in \mc E$, and $\mc L(\theta^{*})$ is a positive semidefinite Laplacian. 
A linearization of the DAE \eqref{eq: closed-loop dynamics -- QS} about the regular set of fixed points $([\theta^*],0)$~and elimination of the algebraic variables gives the reduced Jacobian
\begin{equation*}
	J(\theta^{*}) =
	\underbrace{
 	\begin{bmatrix}
	\!I\! & \boldsymbol{0} \\ \boldsymbol{0} & \!(k \cdot \subscr{D}{tot})^{-1}\!
	\end{bmatrix}
	}_{\triangleq \widetilde Q_{1}}
	\underbrace{
	\begin{bmatrix}
	D_{I}^{-1}\! & \fvec 1 \\
	\fvec 1^{T}\!\!  & \!\subscr{D}{tot}
	\end{bmatrix}
	}_{\triangleq \widetilde Q_{2}}
	\underbrace{
	\begin{bmatrix}
		- \!\subscr{\mc L}{red}(\theta^{*})\! & \boldsymbol{0} \\ \boldsymbol{0} & \!-1\!
	\end{bmatrix}
	}_{\triangleq X}
	\,,
\end{equation*}
where $\subscr{\mc L}{red}(\theta^{*})$ is the Schur complement of $\mc L(\theta^{*})$ with respect to the load entries with indices $\mc V_{I}$. It is known that $\subscr{\mc L}{red}(\theta^{*})$ is again a positive semidefinite Laplacian \cite[Lemma II.1]{FD-FB:11d}. The matrix $\widetilde Q_{1}$ is diagonal and positive definite, and $\widetilde Q_{2}$ is positive semidefinite with $\textup{ker}(\widetilde Q_{2}) = \Span([(D_{I}\fvec 1_{n_{I}}) \,, -1 ])$. 

We proceed via a continuity-type argument. Consider momentarily the perturbed Jacobian $J_{\epsilon}(\theta^{*})$, where $\widetilde Q_{2}$ is replaced by the positive definite matrix
$\widetilde Q_{2,\epsilon}=
\left[\begin{smallmatrix}
\!D_{I}^{-1}\! & \fvec 1 \\
	\!\fvec 1^{T}  & \!\subscr{D}{tot}+ \epsilon \!
\end{smallmatrix}\right]$
with $\epsilon>0$.
The spectrum of $J_{\epsilon}(\theta^{*})$ is obtained from $\widetilde Q_{1} \widetilde Q_{2,\epsilon}  X v = \lambda v$ for some $(\lambda,v) \in \complex \times \complex^{n_I+1}$.
Equivalently, let $y = \widetilde Q_{1}^{-1}v$, then 
\begin{equation*}
	-\, \widetilde Q_{2,\epsilon} \cdot \textup{blkdiag}\!\left( \subscr{\mc L}{red} \,,\, 1/(k \cdot \subscr{D}{tot}) \right) y = \lambda y\,.
\end{equation*}
 The Courant-Fischer Theorem applied to this generalized eigenvalue problem implies that, for $\epsilon > 0$ and modulo rotational symmetry, all eigenvalues $\lambda$ are real and negative. 

Now, consider again the unperturbed case with $\epsilon = 0$. 
{We show that the number of zero eigenvalues for $\epsilon = 0$ equals those for $\epsilon >0$ and thus stability (modulo rotational symmetry) is preserved as $\epsilon \searrow 0$.}
 Recall that $\textup{ker}(\widetilde Q_{2}) = \Span([(D_{I}\fvec 1_{n_{I}}) \,, -1 ])$, and the image of the matrix $\textup{blkdiag}\!\left( \subscr{\mc L}{red} \,,\, 1/(k \cdot \subscr{D}{tot}) \right)$ excludes $\Span([\fvec 1_{n_{I}} \,, 0])$. It follows that $\widetilde Q_{2} \cdot \textup{blkdiag}\!\left( \subscr{\mc L}{red} \,,\, 1/(k \cdot \subscr{D}{tot}) \right) y$ is zero if only if $y \in \textup{Span}([\fvec 1_{n_{I}}\,,0])$ corresponding to the rotational symmetry.
%
We conclude that the number of negative real eigenvalues of $J_{\epsilon}(\theta^{*})$ does not change as $\epsilon \searrow 0$. Hence, the equilibrium $([\theta^*],0)$ of the DAE \eqref{eq: closed-loop dynamics -- QS} is locally exponentially stable. 
\end{IEEEproof}

The CAPI controller \eqref{eq: CAPI},\eqref{eq: CAPI - tuning} preserves the primary power injections while restoring the frequency. 
However, it requires all-to-all communication among the inverters, and a restrictive choice of gains \eqref{eq: CAPI - tuning}. 
To overcome these limitations, we present an alternative controller and a modification of CAPI control. 

\subsection{Distributed Averaging PI (DAPI) Control}
\label{Subsection: Distributed Averaging PI Control}

As third secondary control strategy, consider the {\em distributed averaging proportional integral} (DAPI) controller \cite{JWSP-FD-FB:12u}:
\begin{equation}
	u_{i} = - p_{i} \;,\; k_i \dot p_{i} =D_{i}\dot \theta_{i} + \sum\nolimits_{j \in \mathcal{V}_I} L_{ij}\left(\frac{p_i}{D_i}-\frac{p_j}{D_j}\right)
	.
	\label{eq: DAPI}
\end{equation}
Here, $k_{i}>0$ and $L$ is the Laplacian matrix of a weighted, connected and undirected communication graph between the inverters.
The resulting closed-loop system is then given by
\begin{subequations}\label{eq: secondary control -- SP}
\begin{align}
	0 &= P_i^* - \subscr{P}{e,$i$}(\theta)
	\,,\quad & i \in \mathcal{V}_L  \,,
	\label{eq: secondary control -- SP - 0}
	\\
	D_i\dot{\theta}_i &= {P}_i^* - \subscr{P}{e,$i$}(\theta) - p_{i}
	\,,\quad & i \in \mathcal{V}_I\,,
	\label{eq: secondary control -- SP - 1}
	\\
	 k_i \dot p_{i} &=D_{i}\dot \theta_{i} + \sum\nolimits_{j \in \mathcal{V}_I} L_{ij}\left(\frac{p_i}{D_i}-\frac{p_j}{D_j}\right),\,
	& i \in \mathcal{V}_I \,.
	\label{eq: secondary control -- SP - 2}
\end{align}
\end{subequations}
The following result has been established in earlier work  \cite[Theorem~8]{JWSP-FD-FB:12u} and shows the stability of the closed loop \eqref{eq: secondary control -- SP}.

\begin{theorem}\textbf{(Stability of DAPI-Controlled Network).}
\label{Theorem: Stability of DAPI}
Consider the droop-controlled microgrid \eqref{eq: load flow},\eqref{eq: secondary control} with parameters $P_{i}^{*} \in [0,\overline{P}_i]$, and $D_{i}>0$ for $i \in \mathcal{V}_I$. Let the secondary control inputs be given by \eqref{eq: DAPI} with $k_{i} >0$ for $i \in \mathcal{V}_I$ and a connected communication graph among the inverters $\mc V_{I}$ with Laplacian $L$. The following two statements are equivalent:
\begin{enumerate}
\item[(i)] \textbf{\em Stability of primary droop control:  } the droop control stability condition \eqref{eq: sync condition} holds;

\item[(ii)] \textbf{\em Stability of secondary integral control: } there exists an arc length $\gamma \in [0,\pi/2[$ such that the closed loop \eqref{eq: secondary control -- QS} possesses a locally exponentially stable and unique 
equilibrium manifold $([\theta^*],p^{*}) \subset \overline\Delta_{G}(\gamma)\times \real^{n_{I}}$.
\end{enumerate}
If the equivalent statements (i) and (ii) are true, then $[\theta^{*}]$ {and all injections} are as in Theorem \ref{Thm:Stab}, and $p_{i}^{*} = D_{i} \subscr{\omega}{sync}$, $i \in \mathcal{V}_{I}$. 
\end{theorem}

The DAPI control \eqref{eq: DAPI} regulates the network frequency, requires only a sparse communication network, and preserves the power injections established by primary control. Moreover, the gains $D_{i}>0$ and $k_{i}>0$ can be chosen independently. 

We remark that a higher-order variation of the DAPI control  \eqref{eq: DAPI} (additionally integrating edge flows) can also be derived from a network flow optimization perspective \cite{EM-SHL:14,SY-LC:14}.

%
\subsection{Partial Secondary Control}
\label{Subsection: Partial Secondary Control}

The DAPI and CAPI controller require that {\em all} inverters participate in secondary control. To further reduce the communication complexity and increase the adaptivity  of the microgrid, it is desirable that~only a subset of inverters regulate the frequency. 
To investigate this scenario, we partition the set of inverters as $\mc V_{I} =\mc V_{I_{P}} \cup \mc V_{I_{S}}$, where the action of the $\mathcal{V} _{I_{P}}$ inverters is restricted to primary droop control, and the $\mathcal{V} _{I_{S}}$ inverters perform the secondary DAPI or CAPI control:
\begin{subequations}
\label{eq: partial secondary control}
\begin{align}
	D_i\dot{\theta}_i &= {P}_i^* - \subscr{P}{e,$i$}(\theta)
	\,,\quad & i \in \mathcal{V}_{I_{P}}\,,
	\label{eq: partial secondary control -- 1}
	\\
	D_i\dot{\theta}_i &= {P}_i^*  - \subscr{P}{e,$i$}(\theta)+ u_{i}(t) 
	\,,\quad & i \in \mathcal{V}_{I_{S}}\,,
	\label{eq: partial secondary control -- 2}
\end{align}
\end{subequations}
Observe that the $\mc V_{I_{P}}$ inverters are essentially frequency-dependent loads and the previous analysis applies.
The following result shows that partial secondary control strategies successfully stabilize the microgrid and regulate the frequency.

\begin{theorem}\textbf{(Partially Regulated Networks).}
\label{Theorem: Stability of Partial Secondary Control}
Consider the droop-controlled microgrid with primary and partial secondary control \eqref{eq: load flow},\eqref{eq: partial secondary control} and with parameters $P_{i}^{*} \in [0,\overline{P}_i]$, and $D_{i}>0$ for $i \in \mathcal{V}_I$. For $i \in \mathcal{V}_{I_{S}}$, let the secondary control inputs be given by the CAPI controller \eqref{eq: CAPI}, \eqref{eq: CAPI - tuning} with $|\mc V_{I_{S}}| \geq 1$ and a complete communication graph among the $\mc V_{I_{S}}$ nodes (respectively, by the DAPI controller \eqref{eq: DAPI} with $|\mc V_{I_{S}}| \geq 2$ and a connected communication graph among the $\mc V_{I_{S}}$ nodes).  The following statements are equivalent:
\begin{enumerate}
\item[(i)] \textbf{\em Stability of primary droop control:  } the droop control stability condition \eqref{eq: sync condition} holds;
\item[(ii)] \textbf{\em Stability of partial secondary control: } there is an arc length $\gamma \in [0,\pi/2[$ so that the partially regulated CAPI system \eqref{eq: load flow}, \eqref{eq: CAPI}, \eqref{eq: CAPI - tuning}, \eqref{eq: partial secondary control} (resp.  DAPI system \eqref{eq: load flow}, \eqref{eq: DAPI}, \eqref{eq: partial secondary control})  possesses a locally exponentially stable and unique 
equilibrium manifold $([\theta^*],p^{*}) \subset \overline\Delta_{G}(\gamma)\times \real^{|\mathcal{V}_{I_{S}}|}$.
\end{enumerate}
If the equivalent statements (i) and (ii) hold true, then for $i \in \mc V_{I_{S}} $, the {injections} $\subscr{P}{e,$i$}(\theta^{*})$ are as in Theorem \ref{Thm:Stab} and $p_{i}^{*} = D_{i} \subscr{\omega}{partial}$, where $\subscr{\omega}{partial} = \sum_{i \in \mathcal{V}}P_i^* / (\sum_{i \in \mathcal{V}_{I_S}}D_i)$. For all other inverters $i \in \mc V_{I_{P}}$, we have that $\subscr{P}{e,$i$}(\theta^{*}) = P_{i}^{*}$.
\end{theorem}

\begin{IEEEproof}
The proof for partial CAPI control (respectively, partial DAPI control) is analogous to the proof of Theorem \ref{Theorem: Stability of CAPI} (respectively, \cite[Theorem~8]{JWSP-FD-FB:12u}), while accounting for the partition $\mc V_{I} = \mathcal{V}_{I_P} \cup \mathcal{V}_{I_S}$ in the Jacobian matrices.
\end{IEEEproof}

We now investigate the power sharing properties of partial secondary control. The steady-state  injections at $([\theta^*],p^{*})$ are
\begin{align*}
\subscr{P}{e,$i$}(\theta^*) &= P_{i}^*, \, &  i \in \mc V_{I_{P}} \cup \mc V_{L} \,,
\\
\subscr{P}{e,$i$}(\theta^*) &= P_{i}^* - D_{i}\subscr{\omega}{partial}, \, & i \in \mc V_{I_{S}} \,,
\end{align*}
By applying Theorem \ref{Thm:PowerFlowConstraintsTree}, we obtain the following corollary:

\begin{corollary}\label{Corollary: PowerFlowConstraintsTree}
\textbf{{(Injection Constraints and Power Sharing with Partial
    Regulation).}}  Consider a locally exponentially stable equilibrium 
 $([\theta^*],p^{*}) \subset \overline\Delta_{G}(\gamma)\times
\real^{|\mathcal{V}_{I_{S}}|}$, $\gamma \in {[0,\pi/2[}$, of the partial
    secondary control system \eqref{eq: load flow},\eqref{eq: partial
      secondary control} as in Theorem \ref{Theorem: Stability of Partial
      Secondary Control}. Select the droop coefficients and injection setpoints proportionally. The
    following statements are equivalent:

\begin{enumerate}
\item[(i)] {\textbf{\em Injection constraints: }} $0 \leq \subscr{P}{e,$i$}(\theta^{*}) \leq \overline{P}_i$, $\,\,\,\forall i \in \mathcal{V}_{I_{S}}$;
\item[(ii)] {\textbf{\em Serviceable load: }} $0 \leq - \sum\limits_{j \in \mathcal{V}_{I_P}\cup\mathcal{V}_{L}}P^{*}_{j} \leq \sum\limits_{j \in \mathcal{V}_{I_{S}}}\overline{P}_j\,.$
\end{enumerate}
Moreover, the inverters $\mathcal{V}_{I_{S}}$ performing secondary control share the load proportionally according to their power ratings.
\end{corollary}

{
These results on partial CAPI/DAPI show that only a connected subset of inverters have to participate in secondary control, which further reduces the communication complexity and increases the adaptivity and modularity of the microgrid.}
%
%

\section{Decentralized Tertiary Control Strategies}
\label{Section: Tertiary Control Strategies}

In this section, we examine the tertiary control layer. {In conventional power system operation, the tertiary-level economic dispatch \eqref{eq: AC optimal dispatch} is solved in a centralized way, offline, and with a precise knowledge of the network model and the load profile. In comparison, we show that the economic dispatch \eqref{eq: AC optimal dispatch} is minimized asymptotically by properly scaled droop controllers in a fully decentralized fashion, online, and without a  model.}

{For simplicity and thanks to Lemma~\ref{Lemma: Synchronization Equivalences}, we restrict ourselves to the shifted control system \eqref{eq: primary-controlled system} with the understanding that the optimal asymptotic injections are also obtained by any secondary control that reaches the steady state $u_{i} = -D_{i}\subscr{\omega}{sync}$.}

\subsection{Convex Reformulation of the AC Economic Dispatch}
\label{Subsection: Economic Dispatch}

The main complication in solving the AC economic dispatch optimization
\eqref{eq: AC optimal dispatch} is the nonlinearity and nonconvexity of the
AC injections constraints~\eqref{eq: power flow -- active}. In practical
power system operation, the nonlinear AC injection $\subscr{P}{e}(\theta)$ is often approximated by the
linear {\em DC injection} $\subscr{P}{DC}(\theta)$ with components 
\begin{equation}
	\subscr{P}{DC,$i$}(\theta) = \sum\nolimits_{j=1}^{n}\mathfrak{Im}(Y_{ij}) E_{i} E_{j} (\theta_{i} - \theta_{j})
	\,, \;\;\quad i \in \mathcal{V} 
	\label{eq: power flow -- active -- DC}
	\,.
\end{equation} 
Accordingly, the AC economic dispatch \eqref{eq: AC optimal
  dispatch} is approximated by the corresponding {\em DC economic dispatch} given by
\begin{subequations}
\begin{align}
	\minimize_{\delta \in\real^{n} \,,\, v \in \real^{n_{I}}} 
	& \;\; f(v) = \sum\nolimits_{i \in \mc V_{I}} \frac{1}{2}\alpha_i v_i^2 
	&
	\label{eq: DC optimal dispatch -- cost}\\
	\mbox{subject to} 
	& \qquad P^*_i + v_{i} = \subscr{P}{DC,$i$}(\delta) 
	 &  \forall \; i \in \mc V_{I} \,,
	 \label{eq: DC optimal dispatch -- flow generator}\\
	& \qquad P^*_i = \subscr{P}{DC,$i$}(\delta)   
	 & \forall \;\, i \in \mc V_{L} \,, 
	 \label{eq: DC optimal dispatch -- flow load}\\
	& \qquad | \delta_{i} - \delta_{j} | \leq \supscr{\gamma}{(DC)}_{ij}
	 & \forall \; \{i,j\} \in \mc E \,, 
	 \label{eq: DC optimal dispatch -- security constraint}\\
	 & \qquad \subscr{P}{DC,$i$}(\delta) \in {[0, \overline{P}_i]}
	 & \forall \, i \in \mc V_{I} \,, 
	 \label{eq: DC optimal dispatch -- gen constraint}
\end{align}%
\label{eq: DC optimal dispatch}%
\end{subequations}%
where the DC variables $(\delta,v)$ are distinguished from the AC variables $(\theta,u)$. In formulating the DC economic dispatch \eqref{eq: DC optimal dispatch}, we also changed the line flow parameters from $\supscr{\gamma}{(AC)}_{ij}$ to $\supscr{\gamma}{(DC)}_{ij} \in {[0, \pi/2[}$ for all $\{i,j\} \in \mc E$. The DC  dispatch \eqref{eq: DC optimal dispatch} is a quadratic program with linear constraints and hence convex. 

Typically, the solution $(\delta^{*},v^{*})$ of the DC dispatch \eqref{eq: DC optimal dispatch} serves as proxy for the solution of the non-convex AC dispatch \eqref{eq: AC optimal dispatch}. The following result shows that both problems are equivalent for acyclic networks and appropriate security constraints. 


\begin{theorem}\textbf{(Equivalence of AC and DC Economic Dispatch in Acyclic Networks).}
\label{Theorem: economic dispatch}
Consider the AC economic dispatch \eqref{eq: AC optimal dispatch} and the DC economic dispatch \eqref{eq: DC optimal dispatch} in an acyclic network. The following statements are equivalent:
\begin{enumerate}

	\item \textbf{\em AC feasibility:} the AC economic dispatch problem \eqref{eq: AC optimal dispatch} with parameters $\supscr{\gamma}{(AC)}_{ij} < \pi/2$ for all $\{i,j\} \in \mc E$ is feasible with a global minimizer $(\theta^{*},u^{*}) \in \torus^{n} \times \real^{n_{I}}$; 

	\item \textbf{\em DC feasibility:} the DC economic dispatch problem \eqref{eq: DC optimal dispatch} with parameters $\supscr{\gamma}{(DC)}_{ij} < 1$ for all $\{i,j\} \in \mc E$ is feasible with a global minimizer $(\delta^{*},{v}^{*}) \in \real^{n} \times \real^{n_{I}}$. 
	
\end{enumerate}
If the equivalent statements (i) and (ii) are true, then
$\sin(\supscr{\gamma}{(AC)}_{ij}) = \supscr{\gamma}{(DC)}_{ij}$, $u^{*} =
v^{*}$, $\boldsymbol{\sin}(B^T\theta^*) = B^T\delta^*$, and $f(u^{*}) =
f(v^{*})$ is a global minimum.
\end{theorem}

\begin{IEEEproof}
{The proof relies on the fact that branch flows are unique in an acyclic network: node variables (injections) $P$ can be uniquely mapped to edge variables (flows) $\xi$ via $P = B\xi$.} 

Denote the unique vector of {AC} branch power flows by $\xi = \mc A \sinbf(B^{T} \theta)$; see \eqref{eq: equilibria - primary control - 2}. For an acyclic network, we have  $\mathrm{ker}(B) = \emptyset$, and $\xi \in \real^{n-1}$ can be equivalently rewritten as $\xi = \mc A B^{T} \delta$ for some $\delta \in \real^{n}$. Thus, we obtain 
\begin{equation}
	\mc A \sinbf(B^{T} \theta) = \mc A B^{T}\delta
	\,.
	\label{eq: AC - DC bijection}
\end{equation}
Now, we associate $\delta$ with the angles of the DC flow \eqref{eq: power flow -- active -- DC}, so that \eqref{eq: AC - DC bijection} is a {\em bijective map} between the AC and the DC flows.

Due to the AC security constraints \eqref{eq: AC optimal dispatch -- security constraint}, the sine function is invertible. If the DC security constraints \eqref{eq: DC optimal dispatch -- security constraint} satisfy $\|B^{T}\delta\|_{\infty} \leq \max_{\{i,j\} \in \mc E} \supscr{\gamma}{(DC)}_{ij}<1$, then $B^{T} \theta$ can be uniquely recovered from (and mapped to) $B^{T}\delta$ via \eqref{eq: AC - DC bijection}. Additionally, up to rotational symmetry and modulo $2\pi$, the angle $\theta$ and be uniquely recovered from (and mapped to) $\delta$. Thus, identity \eqref{eq: AC - DC bijection} between the AC and the DC flow serves as a {\em bijective change of variables}  (modulo $2\pi$ and up to rotational symmetry). 

This change of variables maps the AC economic dispatch \eqref{eq: AC optimal dispatch} to the DC economic dispatch \eqref{eq: DC optimal dispatch} as follows.
The AC injections $\subscr{P}{e}(\theta)$ are replaced by the DC injections $\subscr{P}{DC}(\delta)$. The AC security constraint \eqref{eq: AC optimal dispatch -- security constraint} translates uniquely to the DC constraint \eqref{eq: DC optimal dispatch -- security constraint} with $\supscr{\gamma}{(DC)}_{ij} = \sin(\supscr{\gamma}{(AC)}_{ij}) < 1$. 
The AC injection constraint \eqref{eq: AC optimal dispatch -- gen constraint} is mapped to the DC injection constraint \eqref{eq: DC optimal dispatch -- gen constraint}.

Finally, if both problems \eqref{eq: AC optimal dispatch} and \eqref{eq: DC optimal dispatch} are feasible with minimizers $u^{*} = v^{*}$ and $\boldsymbol{\sin}(B^T\theta^*) = B^T\delta^*$, then $f(u^{*}) = f(v^{*})$ is the unique global minimum due to convexity of \eqref{eq: DC optimal dispatch}.
\end{IEEEproof}

Theorem \ref{Theorem: economic dispatch} relies on the bijection \eqref{eq: AC - DC bijection} between  AC and DC flows in acyclic networks \cite{FD-FB:13c,FD-MC-FB:11v-pnas}. For cyclic networks, the two problems \eqref{eq: AC optimal dispatch} and \eqref{eq: DC optimal dispatch} are generally not equivalent, but the DC flow is a well-accepted proxy for the AC~flow.

We now state a rather surprising result: any minimizer of the AC economic dispatch \eqref{eq: AC optimal dispatch} can be achieved by an appropriately designed droop control \eqref{eq: primary control}. Conversely, any steady state of the droop-controlled microgrid \eqref{eq: load flow},\eqref{eq: primary control} is the minimizer of an AC economic dispatch \eqref{eq: AC optimal dispatch} with appropriately chosen parameters. 
The proof {relies on the {\em economic dispatch criterion} \cite{AJW-BFW:96} stating that all marginal costs $\alpha_{i} u_{i}^{*}$ must be identical for the optimal injection}, and it can be extended to the constrained case at the cost of a less explicit relation between the optimization parameters and droop control coefficients.

\begin{theorem}\textbf{(Droop Control \& Economic Dispatch).}
\label{Theorem: Optimal Decentralized Droop Control & Economic Dispatch}
Consider the AC economic dispatch \eqref{eq: AC optimal dispatch} and the shifted control system \eqref{eq: primary-controlled system}.
The following statements are equivalent:
\begin{enumerate}

	\item \textbf{\em Strict feasibility and optimality:} there are parameters $\alpha_{i}>0$, $i \in \mc V_{I}$, and $ \supscr{\gamma}{(AC)}_{ij}<\pi/2$, $\{i,j\} \in \mc E$ such that the AC economic dispatch problem \eqref{eq: AC optimal dispatch} is strictly feasible 
 with global minimizer $(\theta^{*},u^{*}) \in \torus^{n} \times \real^{n_{I}}$. 

	\item \textbf{\em Constrained synchronization:} there exists $\gamma \in {[0, \pi/2[}$ and droop coefficients $D_{i}>0$, $i \in \mc V_{I}$, so that the shifted control system \eqref{eq: primary-controlled system} possesses a unique and locally exponentially stable 
equilibrium manifold $[\theta] \subset \Delta_{G}(\gamma)$ meeting the injection constraints $\subscr{P}{e,$i$}(\theta) \!\in\! {]0, \overline{P}_i[}$,\,$i \in \mc V_{I}$.

\end{enumerate}
If the equivalent statements (i) and (ii) are true, then $[\theta^{*}] = [\theta]$, {$u^{*} = - D \subscr{\omega}{sync}\mathds 1_{n}$}, $\gamma = \max_{\{i,j\} \in \mc E}\supscr{\gamma}{(AC)}_{ij}$, and for some $\beta>0$ 
\begin{equation}
	D_{i} = \beta / \alpha_{i}
	\,, \qquad i \in \mc V_{I} \,.
	\label{eq: optimal droop - strictly feasible}
\end{equation}
\end{theorem}

Theorem \ref{Theorem: Optimal Decentralized Droop Control & Economic Dispatch}, stated for the shifted control system \eqref{eq: primary-controlled system}, can be equivalently stated for the CAPI or DAPI control systems (by Lemma \ref{Lemma: Synchronization Equivalences}). Before proving it, we state a key lemma.%

\begin{lemma}\textbf{(Properties of strictly feasible points).}
\label{Lemma: Properties of feasible points}
If $(\theta^{*},u^{*}) \in \torus^{n} \times \real^{n_{I}}$ is a {strictly feasible} minimizer of the AC economic dispatch \eqref{eq: AC optimal dispatch}, then $u^{*}$ is {\em sign-definite}, that is, all $u_{i}^{*}$, $i \in \mc V_{I}$, have the same sign. Conversely, any strictly feasible pair $(\theta,u) \in \torus^{n} \times \real^{n_{I}}$ of the AC economic dispatch \eqref{eq: AC optimal dispatch} with sign-definite $u$ is {\em inverse optimal} with respect to some $\alpha \in \real_{>0}^{n_{I}}$: there is a set of coefficients $\alpha_{i}>0$, $i \in \mc V_{I}$, such that $(\theta,u)$ is global minimizer of the AC economic dispatch~\eqref{eq: AC optimal dispatch}. 
\end{lemma}

\begin{IEEEproof}
The strictly feasible pairs  of \eqref{eq: AC optimal dispatch} are given by the set of all $(\theta,u) \in \torus^{n} \times \real^{n_{I}}$ satisfying the power flow equations \eqref{eq: AC optimal dispatch -- flow generator}-\eqref{eq: AC optimal dispatch -- flow load} and the strict inequality constraints \eqref{eq: AC optimal dispatch -- security constraint}-\eqref{eq: AC optimal dispatch -- gen constraint}.
Summing all equations \eqref{eq: AC optimal dispatch -- flow generator}-\eqref{eq: AC optimal dispatch -- flow load} yields the necessary solvability condition (power balance constraint) $\sum\nolimits_{i\in\mathcal{V}_I}u_{i} \!=\! -\sum\nolimits_{i\in\mathcal{V}}P_{i}^{*}$.

To establish the necessary and sufficient optimality conditions for \eqref{eq: AC optimal dispatch} in the  strictly feasible case, without loss of generality, we drop the inequality constraints \eqref{eq: AC optimal dispatch -- security constraint}-\eqref{eq: AC optimal dispatch -- gen constraint}. With $\lambda \in \real^{n}$, the Lagrangian $\map{\mc L} {\torus^n \times \real^{n_I} \times \real^n}{\real}$ is given by
\begin{align*}
	\mc{L}(\theta,u,\lambda)
		\!=\!&
		\sum\nolimits _{j \in \mc V_I}\frac{1}{2} \alpha_j u_j^2 + \sum\nolimits_{j \in \mc V_I} \lambda_j  \left(u_j + P^*_j - \subscr{P}{e,$j$}(\theta)\right)
		\\\!&\,+ 
		 \sum\nolimits_{j \in \mc V_L} \lambda_j  \left(P^*_j - \subscr{P}{e,$j$}(\theta) \right) 
		 \,. 
\end{align*}
The necessary KKT conditions \cite{SB-LV:04} for optimality are: 
\begin{subequations} \label{eq: optimality conditions}
\begin{align}
	&\frac{\partial \mc L}{\partial \theta_{i}} = 0\;: \;\;
	0 = \sum\nolimits _{j \in \mc V}\lambda_j \cdot \dfrac{\partial \subscr{P}{e,$j$}(\theta)}{\partial \theta_i}  \,,
	&\forall i \in \mc V \,,\,\,
	\label{eq: optimality conditions -- 1}
	\\
	&\frac{\partial \mc L}{\partial u_{i}} = 0\;: \;\;
	\alpha_i u_i = -\lambda_i \,,
	&\forall i \in \mc V_I \,,
	\label{eq: optimality conditions -- 2}
	\\
	&\frac{\partial \mc L}{\partial \lambda_{i}} = 0\;: \;\;
	-u_i = P_i^* - \subscr{P}{e,$i$}(\theta) \,,
	&\forall i \in \mc V_I \,,
	\label{eq: optimality conditions -- 3}
	\\
	&\frac{\partial \mc L}{\partial \lambda_{i}} = 0\;: \;\;
	0 = P_i^* - \subscr{P}{e,$i$}(\theta) \,,
	&\forall i \in \mc V_L \,.
	\label{eq: optimality conditions -- 4}
\end{align}
\end{subequations}
Since the AC economic dispatch \eqref{eq: AC optimal dispatch} is equivalent to the convex DC dispatch (see Theorem \ref{Theorem: economic dispatch}), the KKT conditions \eqref{eq: optimality conditions} are also sufficient for optimality.
In vector form, \eqref{eq: optimality conditions -- 1} reads as 
$\fvec 0_{n} = \lambda^{\transpose}\partial \subscr{P}{e}(\theta)/\partial \theta$, where the {\em load flow Jacobian} is given by symmetric Laplacian
$\partial \subscr{P}{e}(\theta)/\partial \theta = B \diag(\{a_{ij}\}_{\{i,j\} \in \mc E})B^{T}$
with strictly positive weights $a_{ij} = \mathfrak{Im}(Y_{ij}) E_{i} E_{j} \cos(\theta_{i}-\theta_{j})$ (due to strict feasibility of the security constraint \eqref{eq: AC optimal dispatch -- security constraint}). 

It follows that  $\lambda \in \mathds 1_{n}$, that is, $\lambda_i = \widetilde{\lambda} \in \real$ for all $i \in \mc V$ and for some $\widetilde{\lambda} \in \real$. 
Hence, condition \eqref{eq: optimality conditions -- 2} reads as the economic dispatch criterion (identical marginal costs) $u_{i} = -\widetilde{\lambda}/\alpha_{i}$ for all $i \in \mc V_I$, and the conditions \eqref{eq: optimality conditions -- 3}-\eqref{eq: optimality conditions -- 4} reduce to
\begin{subequations} \label{eq: equivalent optimality conditions -- rewritten}
\begin{align}
\widetilde{\lambda}/\alpha_{i} &= P_i^* - \subscr{P}{e,$i$}(\theta) \,, &\forall i \in \mc V_I \,,
\label{eq: equivalent optimality conditions -- rewritten -- 1}
\\
0 &= P_i^* - \subscr{P}{e,$i$}(\theta) \,, &\forall i \in \mc V_L \,.
\label{eq: equivalent optimality conditions -- rewritten -- 2}
\end{align}
\end{subequations}
By summing all equations \eqref{eq: equivalent optimality conditions -- rewritten}, we obtain the constant $\widetilde \lambda$ 
as $	\widetilde \lambda = {\sum\nolimits_{i\in\mathcal{V}}P_{i}^{*}}/{\sum\nolimits_{i\in\mathcal{V_{I}}}\alpha_{i}^{-1}}$.
The minimizers are $u_{i}^{*} = -\widetilde{\lambda}/\alpha_{i}$ and $\theta^{*}$  determined from \eqref{eq: equivalent optimality conditions -- rewritten}. 
It follows that $u^{*}$ is sign-definite. 

By comparing the (strict) optimality conditions \eqref{eq: equivalent optimality conditions -- rewritten} with the (strict) feasibility conditions \eqref{eq: AC optimal dispatch -- flow generator}-\eqref{eq: AC optimal dispatch -- flow load}, 
it follows that any strictly feasible pair $(\theta,u)$ with sign-definite $u$ is inverse optimal for the coefficients $\alpha_{i} = -\beta/u_{i}$ with some $\beta>0$.
\end{IEEEproof}

\begin{IEEEproof}[Proof of Theorem \ref{Theorem: Optimal Decentralized Droop Control & Economic Dispatch}]%
\textbf{(i) $\implies$ (ii):}
 If the AC economic dispatch~\eqref{eq: AC optimal dispatch} is strictly feasible, then its minimizer $(\theta^{*},u^{*})$ is global (Theorem \ref{Theorem: economic dispatch}), and the optimal inverter injections are $	\supscr{P}{opt}_{i} = \subscr{P}{e,$i$}(\theta^{*}) = P_{i}^{*} + u_{i}^{*}$ with sign-definite $u^{*}$ (Lemma \ref{Lemma: Properties of feasible points}). 
 Since the power flow equations \eqref{eq: AC optimal dispatch -- flow generator}-\eqref{eq: AC optimal dispatch -- flow load} and the strict inequality constraints \eqref{eq: AC optimal dispatch -- security constraint}-\eqref{eq: AC optimal dispatch -- gen constraint} are met, $\supscr{P}{opt}_{i} \in {]0,\overline{P}_{i}[}$, $[\theta^{*}] \subset \Delta_{G}(\gamma)$ with $\gamma = \max_{\{i,j\} \in \mc E}\supscr{\gamma}{(AC)}_{ij}\!$, and the vector of load and source injections $(P_L^*,\supscr{P}{opt}_{I})$ is a $\gamma$-feasible injection setpoint. 

 By Theorem \ref{Theorem: power injection setpoint design} and identity \eqref{eq: DroopSetpoint}, the droop coefficients
$
	D_i = -\beta ({P}_i^{*} - \supscr{P}{opt}_{i}) = \beta u_{i}^{*} 
$, $i \in \mathcal{V}_I$, 
guarantee that the shifted control system \eqref{eq: primary-controlled system} possesses an equilibrium manifold $[\theta]$ satisfying $P_{\rm e}(\theta) = P^{\rm opt} = P_{\rm e}(\theta^*)$. For $\beta u_{i}^{*}>0$ (recall $u^{*}$ is sign-definite), $[\theta]$ is locally exponentially stable by Theorem \ref{Thm:Stab}. Finally, the identity $P_{\rm e}(\theta) = P_{\rm e}(\theta^*)$ shows that $[\theta^{*}] = [\theta]$.

\textbf{ (ii) $\implies$ (i):}
Any equilibrium manifold $[\theta] \subset \Delta_{G}(\gamma)$ as in (ii) is a $\gamma$-feasible power injection setpoint with 
\begin{subequations}
\begin{align}
	 &\widetilde P_{i} = P^*_i - D_{i} \subscr{\omega}{sync} = \subscr{P}{e,$i$}(\theta) 
	 &  \forall \; i \in \mc V_{I} \,,
	 \label{eq: steady state inverse optimality - 1}
	 \\
	&\widetilde P_{i} = P^*_i= \subscr{P}{e,$i$}(\theta)   
	 & \forall \;\, i \in \mc V_{L} \,, 
	 \label{eq: steady state inverse optimality - 2}
	\\
	& | \theta_{i} - \theta_{j} | < \gamma
	 & \forall \; i,j \in \mc E \,, 
	\\
	&\subscr{P}{e,$i$}(\theta) \in {]0, \overline{P}_i[}
	 & \forall \, i \in \mc V_{I} \,.
\end{align}%
\end{subequations}%
Hence, any $\theta \in [\theta]$ is strictly feasible for the economic dispatch \eqref{eq: AC optimal dispatch} if we identify $\theta^{*}$ with $\theta$ (modulo symmetry), $\gamma$ with $\max_{\{i,j\} \in \mc E}\supscr{\gamma}{(AC)}_{ij}$, and $u_i^{*}$ with $-D_{i} \subscr{\omega}{sync}$ (modulo\,scaling). %
Since $u_i^{*}$ is sign-definite, the claim follows from Lemma~\ref{Lemma: Properties of feasible points}. 

In the strictly feasible case, a comparison of the stationarity conditions \eqref{eq: steady state inverse optimality - 1}-\eqref{eq: steady state inverse optimality - 2} and the optimality conditions \eqref{eq: equivalent optimality conditions -- rewritten}~gives 
$
D_{i} \subscr{\omega}{sync} = -u_{i}^{*} = \widetilde{\lambda}/\alpha_{i}, 
$
where $\subscr{\omega}{sync}$ and $\widetilde{\lambda}$ are constant. Since the droop gains are defined up to scaling, we obtain \eqref{eq: optimal droop - strictly feasible}. 
\end{IEEEproof}

\begin{remark}[\bf Selection of droop coefficients]
The equivalence revealed in Theorem \ref{Theorem: Optimal Decentralized Droop Control & Economic Dispatch} suggests the following guidelines to select the droop coefficients: large coefficients $D_{i}$ for desirable (e.g., economic or low emission) sources with small {cost coefficients}~$\alpha_{i}$\,; and vice versa. These insights can also be connected to the  proportional power sharing objective: if each $P_{i}^{*}$ and $1/\alpha_{i}$ are selected proportional to the rating $\overline{P}_{i}$, that is, $\alpha_{i}{\overline P_{i}} = \alpha_{j}{\overline P_{j}}$ and $P_i^*/\overline{P}_i = P_j^*/\overline{P}_j$, then the associated droop coefficients \eqref{eq: optimal droop - strictly feasible} equal those in \eqref{eq: prop choice of droop coeff} for load sharing. 
\oprocend
\end{remark}

{\begin{remark}[\bf Beyond quadratic objectives and linear droop slopes]
\label{Remark: Beyond quadratic objective functions}
As shown in Theorem~\ref{Theorem: Optimal Decentralized Droop Control & Economic Dispatch}, the steady states of a microgrid \eqref{eq: load flow} with linear droop control \eqref{eq: primary control} are related one-to-one to the global optimizers of the economic dispatch \eqref{eq: AC optimal dispatch} with quadratic objective  \eqref{eq: AC optimal dispatch -- cost}. If analogous proof methods are carried out for a general objective function $f(u) = \sum\nolimits_{i \in \mc V_{I}} C_{i}(u_{i})$ with strictly convex and continuously differentiable functions $C_{i}$, the associated optimal droop controllers \eqref{eq: primary control} need to have nonlinear frequency-dependent droop slopes given by $D_{i} = {(C_{i}^{\prime})}^{-1}(\dot\theta_{i})$. Conversely, practically employed droop curves with frequency deadbands and saturation \cite[Chapter 9]{JM-JWB-JRB:08} can be related to non-smooth and barrier-type costs \cite{SY-LC:14}.

From such an optimization perspective, the primary dynamics \eqref{eq: primary-controlled system} are a primal algorithm converging to a steady-state satisfying the  optimality conditions \eqref{eq: optimality conditions}. 
Likewise, second-order or integral control dynamics {can be interpreted} as primal-dual algorithms, as shown for related systems 
\cite{MA-DVD-KHJ-HS:13,MA-DVD-HS-KHJ:14,EM-SHL:14,SY-LC:14,XZ-AP:13,CZ-UT-NL-SL:13,NL-CL-ZC-SHL:13}.
\oprocend
\end{remark}
}

%
%
%
%

\section{Simulation Case Study}
\label{Section: Simulation Scenario}

{
Throughout the past sections we demonstrated that the CAPI control \eqref{eq: secondary control}, \eqref{eq: CAPI} and DAPI control \eqref{eq: secondary control}, \eqref{eq: DAPI} with properly scaled coefficients can simultaneously address primary, secondary, and tertiary-level objectives in a plug-and-play fashion. Both strategies rely on simple distributed and averaging-based PI controllers that do not require a hierarchical implementation with time-scale separations and detailed system knowledge.

We illustrate the performance of our controllers via simulation of the IEEE 37 distribution grid \cite{IEEE37:10} shown in Fig.~\ref{Fig:Sim}. 
After an islanding event, the distribution grid is disconnected from the transmission network, and distributed generators must ensure stability while~regulating the frequency and sharing the demand. 
The communication network among the distributed generators is shown in dotted blue. Of the 16 sources, 8 have identical power ratings, while the remaining 8 are rated for twice as much power. 
To demonstrate the robustness of our controllers beyond our theoretical results, we use the coupled and lossy power flow \cite{PK:94} in place of the lossless\,and decoupled equations \eqref{eq: power flow}. On the reactive power side we control the inverter voltages via the \emph{quadratic voltage-droop controller}~\cite{JWSP-FD-FB:13h}
\begin{equation}
\tau_i\dot{E}_i = -C_iE_i(E_i-E_i^*) - Q_{\mathrm{e},i},\,\quad i \in \mathcal{V}_I \,,
\label{eq: quadratic droop}
\end{equation}
where $E_i^* > 0$ is the nominal voltage, $C_i>0$ and $\tau_i > 0$ are gains, and $Q_{\mathrm{e},i} \in \real$ is the reactive power injection \eqref{eq: power flow -- reactive}.

We compare the performance of primary droop control \eqref{eq: primary control}, decentralized secondary integral control \eqref{eq: decentralized control} at every source, and the DAPI control \eqref{eq: DAPI} after a step change at a single load. 
We choose as objective proportional power sharing \eqref{eq: prop power sharing} or equivalently economic dispatch \eqref{eq: AC optimal dispatch} with coefficients $\alpha_{i} = 1/\overline{P_{i}}$, and the droop control coefficients are obtained accordingly from \eqref{eq: optimal droop - strictly feasible}. 
The secondary controller time constants have been randomized to model a true plug-and-play scenario, where only  communication channels have been established without any tuning of control gains. 
As can be seen in Fig.~\ref{Fig:Freqs}, primary droop control gives rise to a frequency deviation,  while both decentralized integral control and DAPI control quickly regulate the frequency with similar closed-loop dynamics, but drastically different power injections (Fig.~\ref{Fig:Powers}) and marginal costs (Fig.~\ref{Fig:MCosts}). All three controllers give rise to very similar voltage dynamics (Fig. \ref{Fig:Fb} depicts the case for DAPI control) together with the quadratic droop controller \eqref{eq: quadratic droop} on the reactive power side.
 In comparison to decentralized integral control (Fig. \ref{Fig:Pa} and \ref{Fig:Ma}), DAPI control ensures proportional power sharing (Fig. \ref{Fig:Pb}) and economic optimality, as seen from the asymptotically equal marginal generation costs (Fig.~\ref{Fig:Mb}).

\begin{figure}[!t]
	\centering
	\begin{subfigure}[!ht]{0.5\columnwidth}
		\includegraphics[width=\columnwidth]{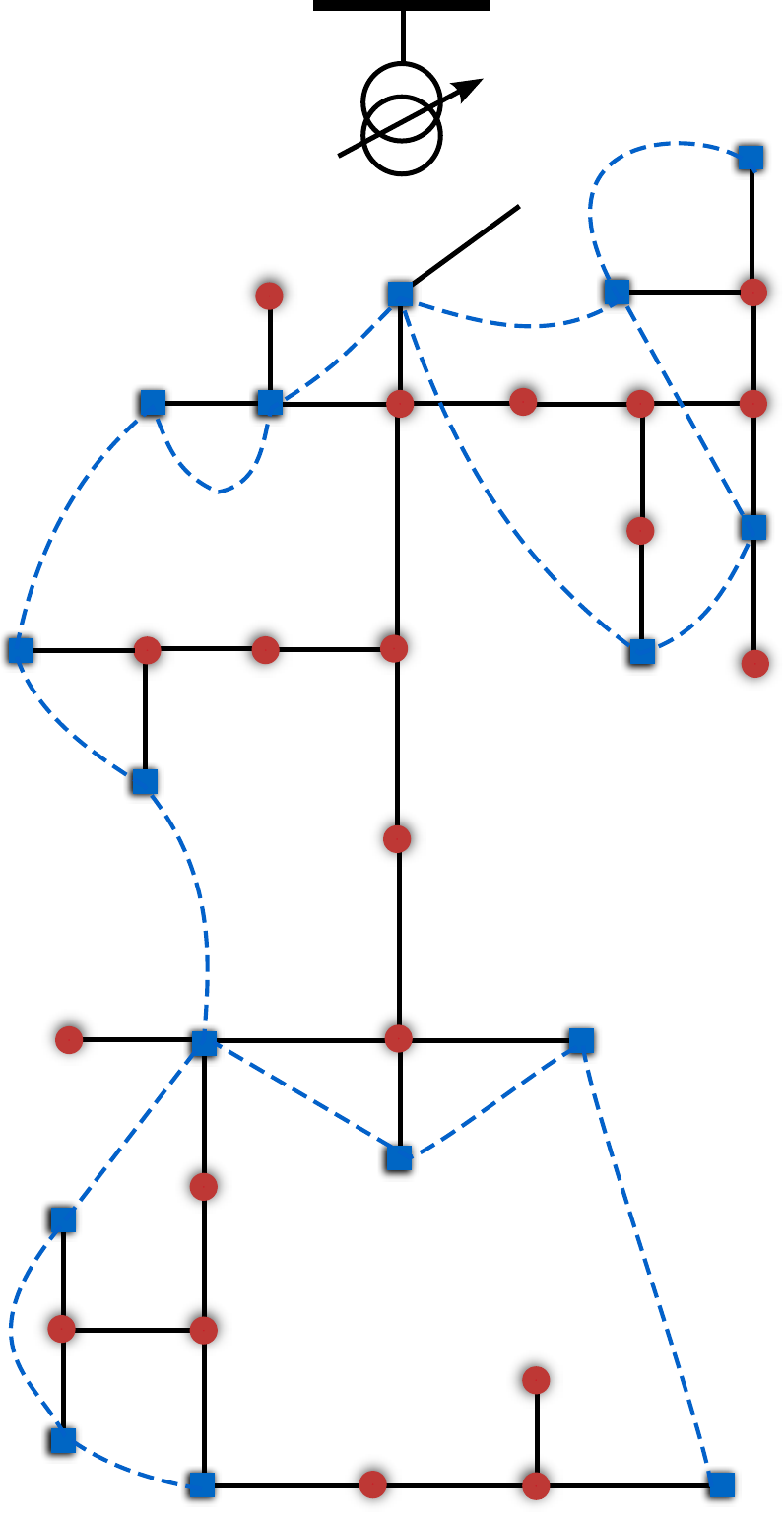}
		\caption{Islanded IEEE37 feeder}
		\label{Fig:Sim}
	\end{subfigure}~
	\begin{subfigure}[!ht]{0.5\columnwidth}
		\begin{subfigure}[!ht]{\columnwidth}
			\includegraphics[width=\columnwidth]{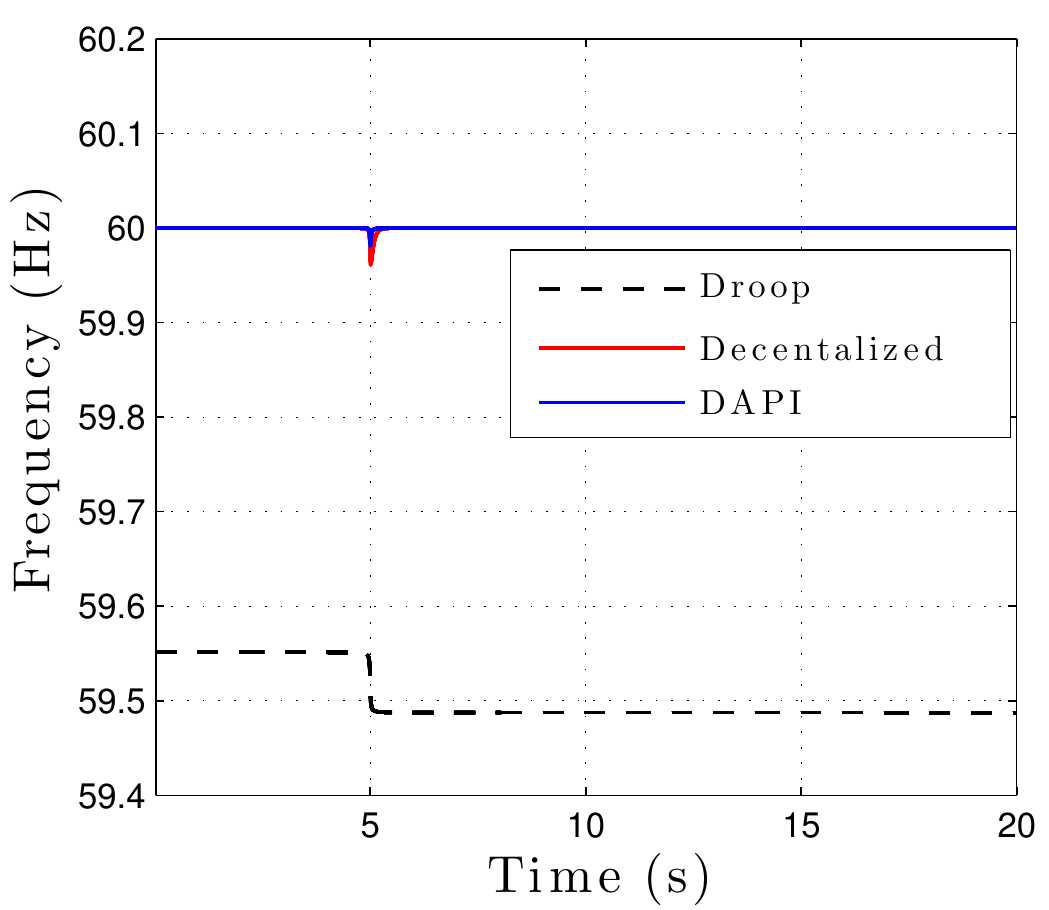}
			\caption{Frequency dynamics}
			\label{Fig:Fa}
		\end{subfigure}
		\begin{subfigure}[!ht]{\columnwidth}
			\includegraphics[width=\columnwidth]{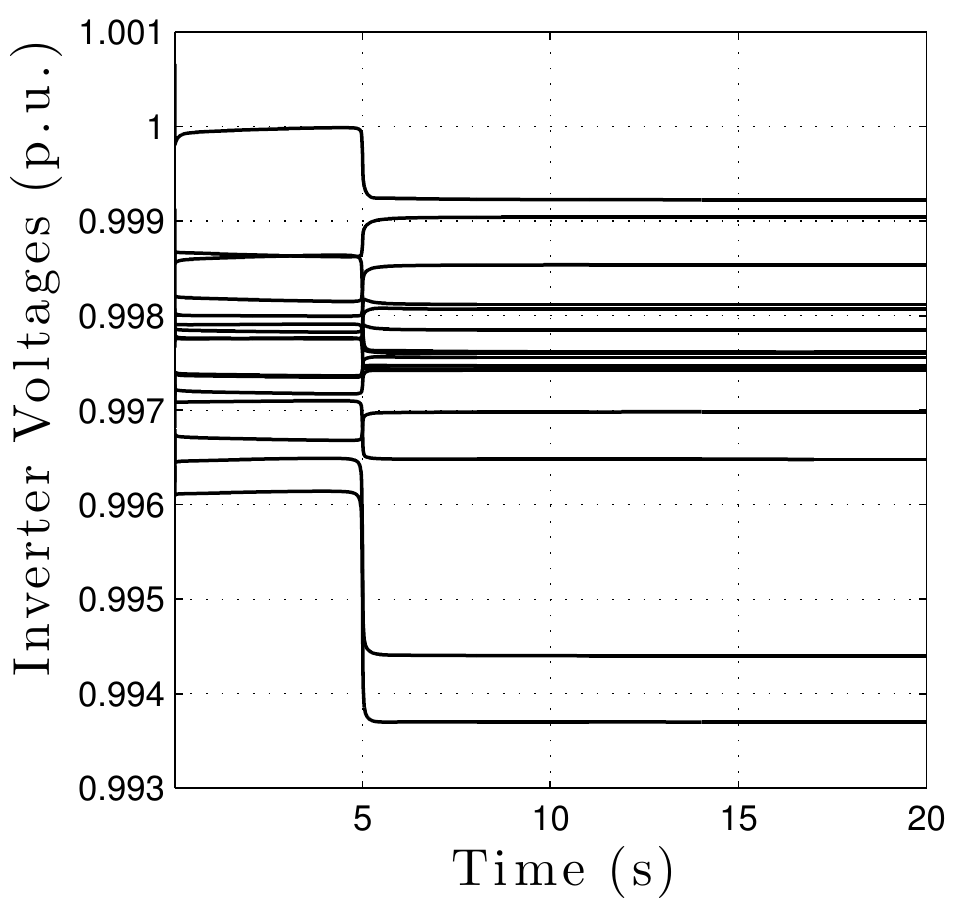}
			\caption{Voltage magnitude dynamics}
			\label{Fig:Fb}
		\end{subfigure}
	\end{subfigure}
	\caption{Depiction of the islanded IEEE 37 microgrid with loads (red nodes) and generation units (blue nodes) interfaced with droop-controlled inverters; and comparison of frequency and voltage regulation under primary droop control  \eqref{eq: primary control}, decentralized integral control \eqref{eq: decentralized control}, and DAPI \eqref{eq: DAPI}  frequency control.}\label{Fig:Freqs}	
\end{figure}
\begin{figure}[!h]
        \centering
        \begin{subfigure}[!ht]{0.45\columnwidth}
                \includegraphics[width=\columnwidth]{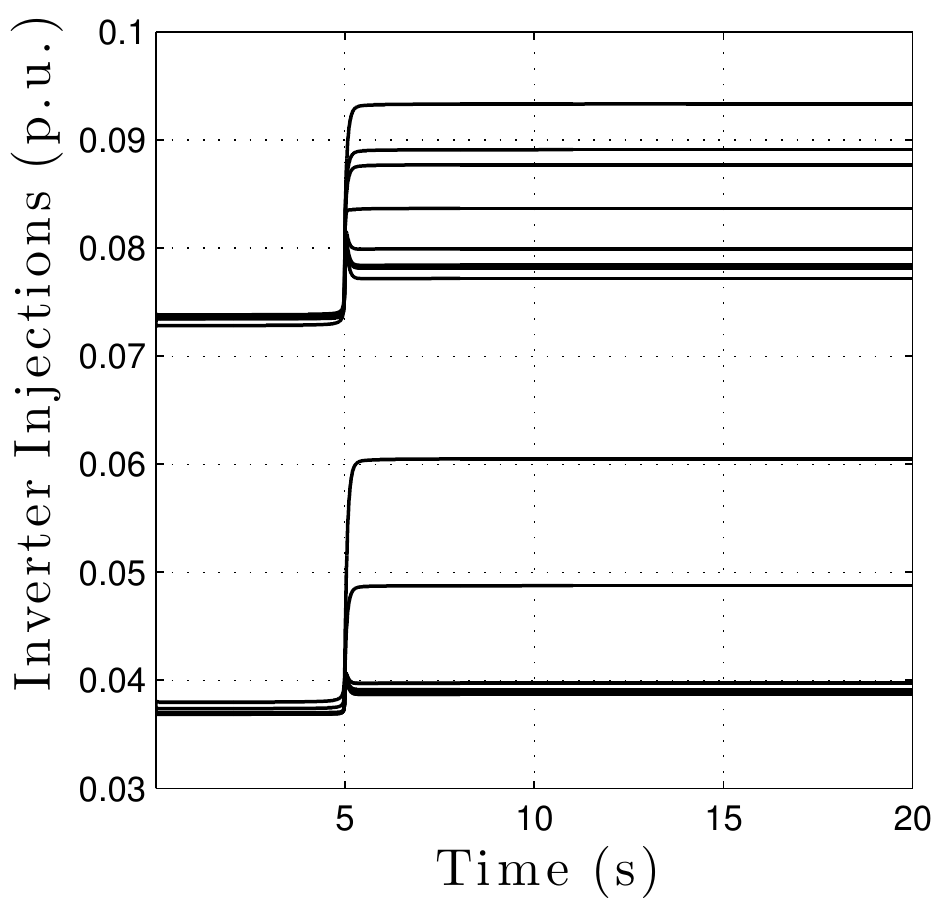}
                \caption{Decentralized control}
                \label{Fig:Pa}
        \end{subfigure}~
        \begin{subfigure}[!ht]{0.45\columnwidth}
                \includegraphics[width=\columnwidth]{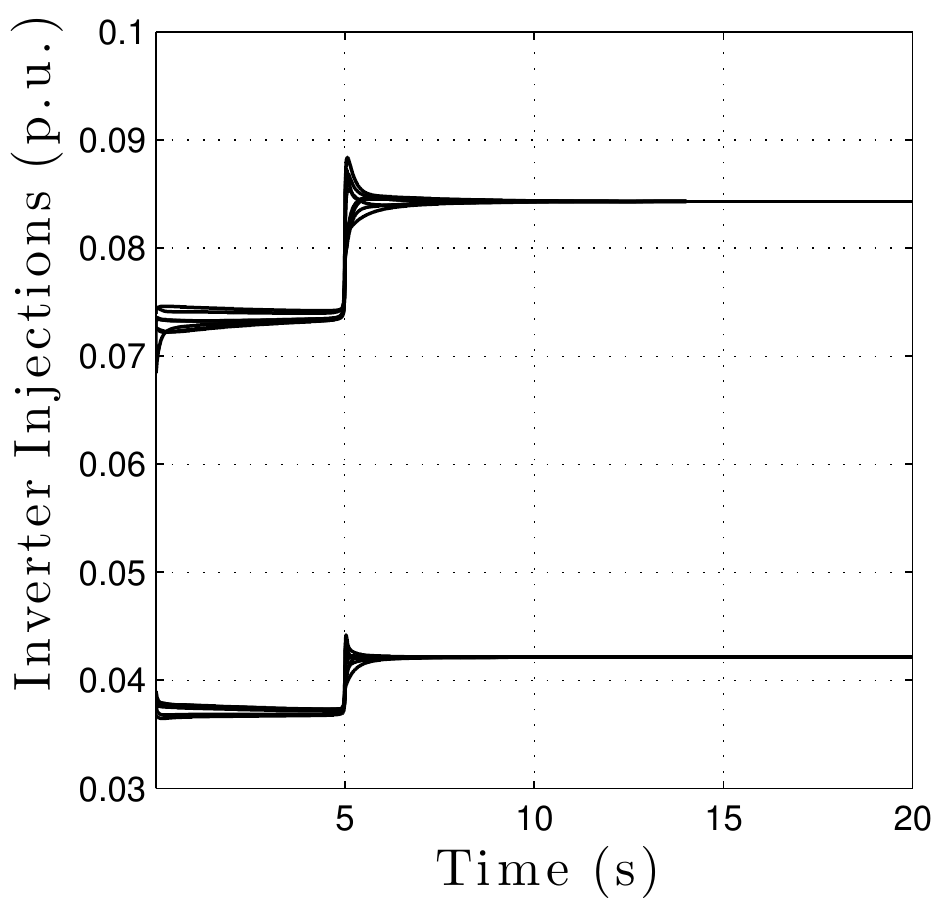}
                \caption{DAPI control}
                \label{Fig:Pb}
		\end{subfigure}
        \caption{Dynamics of power injections $P_{\mathrm{e},i}$ after a step change in load.}\label{Fig:Powers}
\end{figure}
\begin{figure}[!h]
        \centering
        \begin{subfigure}[!ht]{0.45\columnwidth}
                \includegraphics[width=\columnwidth]{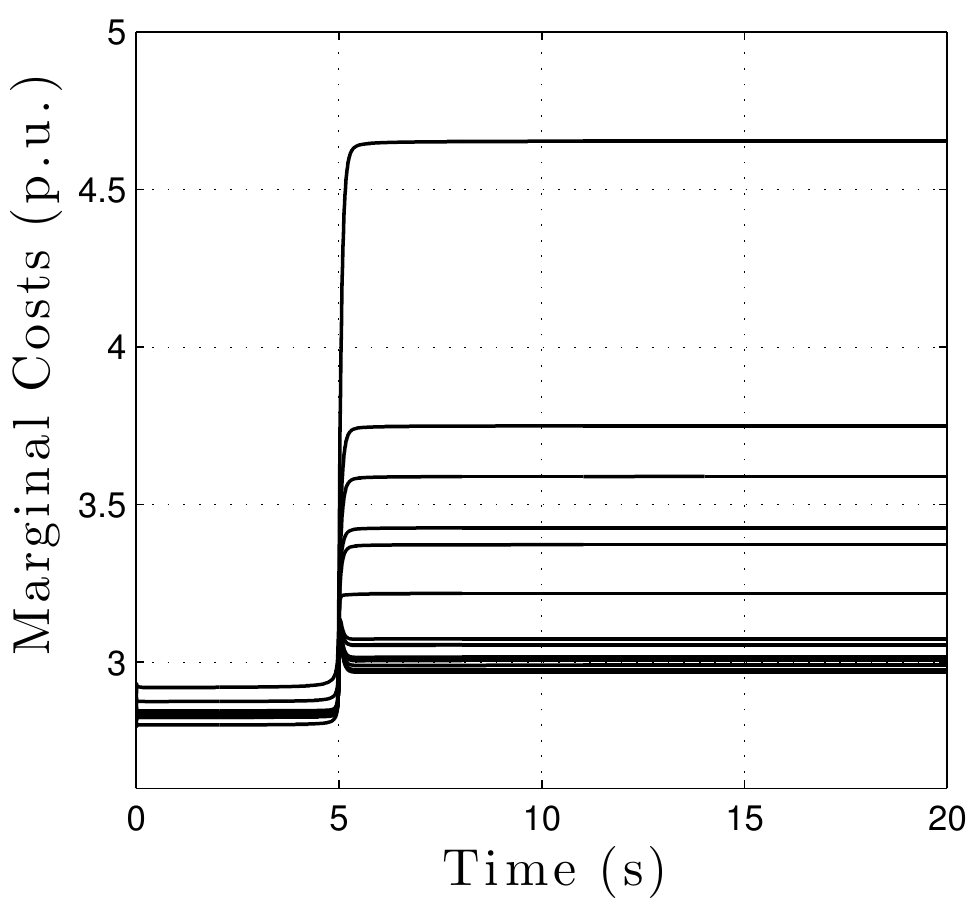}
                \caption{Decentralized control}
                \label{Fig:Ma}
        \end{subfigure}~
        \begin{subfigure}[!ht]{0.45\columnwidth}
                \includegraphics[width=\columnwidth]{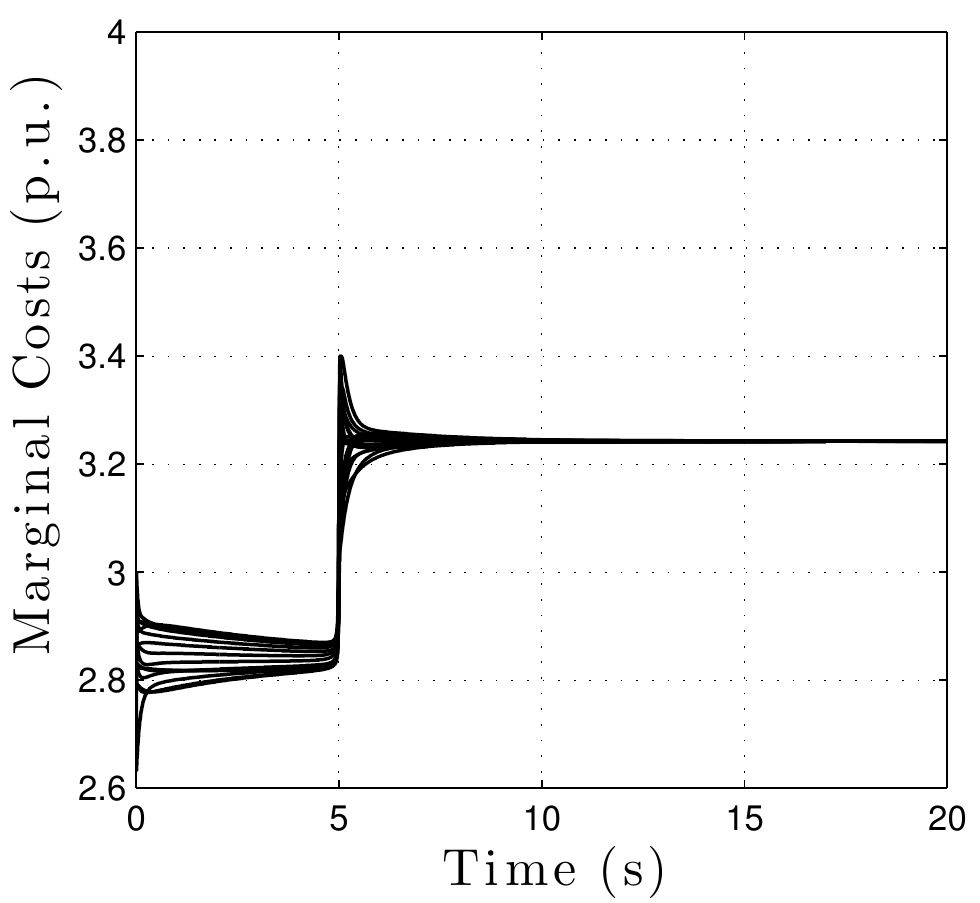}
                \caption{DAPI control}
                \label{Fig:Mb}
		\end{subfigure}
        \caption{Dynamics of marginal costs $\alpha_iu_i$ after a step change in load.}\label{Fig:MCosts}
\end{figure}

}

%
\section{Conclusions}
\label{Section: Conclusions}

We studied decentralized and distributed primary, secondary, and tertiary control strategies in microgrids and illuminated some connections between them. Thereby, we relaxed some restrictions regarding the information structure and time-scale separation of conventional hierarchical control strategies adapted from transmission-level networks to make them more applicable to microgrids and distribution-level applications.

While this work is a first step towards an understanding of the
interdependent control loops in hierarchical microgrids, several
complicating factors have not been taken into account. In particular, our
analysis is {only local and so far} formally restricted to acyclic networks with {constant resistance-to-reactance ratios}. Moreover, future work needs to consider more
detailed models including reactive power flows, voltage dynamics, and
ramping constraints on the inverter injections. {In preliminary work \cite{CZ-EM-FD:15} we extend the present analysis to cyclic networks possibly with higher-order generator dynamics in transmission grid settings, and we provide some first guarantees on the region of attraction. Finally, another interesting direction for future work is to remove the idealistic communication assumptions and resort to sampled or event-triggered schemes in presence of delays. Event-triggered or deadband-enforcing control could also be useful for relaxing frequency regulation by ignoring sufficiently small deviations.}


\section*{Acknowledgments}

The authors wish to thank H. Bouattour, J. M. Guerrero, Q.-C. Zhong, A. Dom\'inguez-Garc\'ia, N. Ainsworth, and M. Andreasson for insightful discussions and sharing their preprints.}


\renewcommand{\baselinestretch}{0.945}
\bibliographystyle{IEEEtran}
\bibliography{alias,Main,FB,New}

\end{document}